\documentclass[11pt]{amsart}
\usepackage[letterpaper, margin=1in]{geometry}
\usepackage{graphicx}
\everymath{\displaystyle} 
\usepackage{graphicx}

\usepackage{tikz,amsfonts,amsmath,amssymb,mathrsfs}
\usepackage{enumitem}

\newtheorem{thm}{Theorem}[section]%

\newtheorem{dfn}[thm]{Definition}
\newtheorem{pro}[thm]{Proposition}
\newtheorem{cor}[thm]{Corollary}
\newtheorem{lemme}[thm]{Lemma}
\newtheorem{rmk}[thm]{Remark}

\renewcommand{\phi}{\varphi}

\numberwithin{equation}{section}

\newcommand{\di}{\operatorname{div}}

\newcommand{\lebd}{\mathcal{L}^d}

\newcommand \p {\subset  }
\newcommand {\bd}{\begin{dfn}}
\newcommand \ed{\end{dfn}}
\newcommand \be{\begin{equation}}
\newcommand \ee{\end{equation}}
\newcommand{\intom}{\int_\Omega }
\newcommand{\cqfd}{\begin{flushright}$\square$\end{flushright}}
\newcommand \bt{\begin{thm}}
\newcommand \et{\end{thm}}
\newcommand \bp{\begin{pro}}
\newcommand \ep{\end{pro}}
\newcommand \bc{\begin{cor}}
\newcommand \ec{\end{cor}}
\newcommand \bl{\begin{lemme}}
\newcommand \El{\end{lemme}}
\newcommand \br{\begin{rmk}}
\newcommand \Er{\end{rmk}} 
\newcommand{\s}{^{-1}}
\newcommand \R{\mathbb{R}  }

\newcommand \Rd{\mathbb{R}^d  }

\newcommand \N{\mathbb{N}  }

 \newcommand{\mat}[2]{\R^{#1\times #2}}

\newcommand \seq[3] { \{ #1_{#2}\}^{\infty } _{#2=#3} }

\newcommand \Lii[1]{ \lim _{#1\rightarrow \infty} }

\newcommand{\w}[3]{W^{#1,#2}(#3)}
\newcommand{\wo}[3]{W_0^{#1,#2}(#3)}
\newcommand{\I}{\infty}
 
\newcommand{\el}[2]{  L^{#1}(#2)  }

\newcommand{\pd}{\partial}
\newcommand{\pair}[2]{\left\langle{#1},{#2}\right\rangle}

\newcommand{\Set}{\mathcal{S}}

\begin{document}

\title[On the Uniqueness of minimizers for a class of variational problems]{ On the uniqueness of minimizers for a class of variational problems with  polyconvex integrand.}

\author{Romeo Awi}
\address{Department of Mathematics, Hampton University, Hampton, Virginia 23668}
\email{romeo.awi@hamptonu.edu}
\thanks{Most of the work presented in this paper was carried out while R.A. was a postdoctoral fellow at the Institute for Mathematics and its Applications during the IMA's annual program on Control Theory and its Applications.
}

\author{Marc Sedjro}
\address[Marc Sedjro]{
        African Institute for Mathematical Sciences(AIMS), Tanzania, \\
 Plot 288 No 288, Makwahiya Street, Regent Estate, P.O. Box 106077, Dar Es Salaam, Tanzania.}
\email{sedjro@aims.ac.tz}
\thanks{M.S. gratefully acknowledges the support of the King Abdullah University of Science and Technology.}


\begin{abstract}
We prove existence and uniqueness of minimizers for a family of energy
functionals that arises in Elasticity and involves polyconvex
integrands over a certain subset of displacement maps. This work
extends previous results by Awi and Gangbo to a larger class of
integrands. We are interested in Lagrangians of the form $
L(A,u)=f(A)+H(\det A)-F\cdot u $. Here the strict convexity condition
on $ f $ and $ H $ have been relaxed merely to a convexity condition.
Meanwhile, we have allowed the map $ F $ to be non-degenerate.
First, we study these variational problems over displacements for which
the determinant is positive. Second, we consider a limit case in which
the functionals are degenerate. In that case, the set of admissible
displacements reduces to that of incompressible displacements which are
measure preserving maps. Finally, we establish that the minimizer over
the set of incompressible maps may be obtained as a limit of minimizers
corresponding to a sequence of minimization problems over general
displacements provided we have enough regularity on the dual problems.
We point out that these results do not rely on the direct methods of the calculus
of variations.
\end{abstract}

\maketitle

\section{Introduction}
We are interested in Euler-Lagrange equations,  existence and uniqueness of minimizers  for some problems in the vectorial calculus of variations emanating from elasticity theory. These variational problems are related  to an open problem in Partial Differential Equations that we describe as follows : let $T>0$ and let   $\Omega$ and $\Lambda$ be two open subsets of  $\mathbb{R}^d$; suppose  that $\mathbf{u}_0$ is a diffeomorphism  between $\Omega $ and $\Lambda$;  we seek $\mathbf{u} :\Omega\times (0,T)\longrightarrow \mathbb{R}^d $ such that  $\mathbf{u}(\cdot,t) (\Omega)=\Lambda$ for each $t$ and   
\begin{equation} \label{eqn: PDE}
 \begin{cases}
  \mathbf{u}_t=\text{div}_x D_\xi L(\nabla \mathbf{u})&\text{ on }\Omega\times (0,T),
  \\
 \mathbf{u}(0, \cdot)=\mathbf{u}_0  &\text{ on }\Omega,
 \end{cases}
\end{equation}
 in the sense of distributions. In \eqref{eqn: PDE}, we assume that the map  $\R^{d\times d}\ni\xi\mapsto L(\xi)$ is quasiconvex. We refer the reader to \cite{awi-gan}, %
 \cite{egs}, %
\cite {Evans20120339}, %
\cite {Gan_uecsc}, %
 and %
\cite{gan-wes} for further details  on these gradient flows.
Understanding variational problems associated to  the time-discretization of (\ref{eqn: PDE}) is arguably an important step toward the construction of a solution. In that regard,  several partial results are available in the literature  (See for instance
\cite{egs} and  %
\cite {Evans20120339}).%

In \cite {awi-gan}, %
the authors have focused on  a class of Lagrangians  that arises  in elastic materials. More precisely, they have considered polyconvex Lagrangians of the form $\xi\mapsto L(\xi)=f(\xi)+H(\det\xi)$.  Here $f$ is a $C^1(\mathbb{R}^d  )$ strictly convex function with $p$-th order growth, and the map $H$ is a $C^1(0,\infty)$ convex  function that satisfies 
\begin{equation} 
\label{eq:cond}
  \lim _{t\rightarrow 0^+}  H(t)= \lim _{t\rightarrow \infty} \frac{H(t)}{t}=+\infty.
\end{equation}
As a result, a variational problem  emerges from the time discretization  and has a relaxation that takes the general form :
\begin{equation} 
\label{eq:pro1}
\min \left\{ \int_\Omega  \left( f(\nabla u)+H(\beta)-F\cdot u \right)dx ;\; (u,\beta)\in \mathscr{U}\right\}
\end{equation}
where $F\in\el 1{\Omega,\R^{d}}$ and
\begin{equation}
\begin{split}
\mathscr{U}=\Big\{(u,\beta): & u\in W^{1,p}(\Omega,\bar\Lambda) ,\; \beta:\Omega\to [0,\I);
\\
&  \int_\Omega  l(u)\beta dx=\int_\Lambda l(y)dy;\forall l\in C_c(\mathbb{R}^d  ) \Big\}.
\end{split}
\end{equation}
 Although the  existence of minimizers in \eqref{eq:pro1}  follows from the direct methods in the calculus of variations, the uniqueness is a rather challenging problem. Indeed, because of \eqref{eq:cond} and the non-convexity of the integrand, standard techniques in calculus of variations do not apply. 

To bypass these difficulties, the authors of  \cite {awi-gan} %
have introduced a pseudo-pro\-jected gradient operator $\mathscr U_\Set\ni u\mapsto \nabla_\Set u$ defined as follows : for a given $u\in \mathscr U_\Set$, the map $\nabla_\Set u$ is the unique minimizer of
\[
 \intom f(G)dx
\]
over 
\[
\mathcal G_\Set(u):=
\left\{
G\in \el p {\Omega,\R^{d\times d}}:
\intom u\di \phi\; dx =-\intom \pair{G}{\phi}dx\;
\forall \phi\in \Set  
\right\}.
\]
Here, $\Set$ is a finite-dimensional subspace of $ W_0^{1,q}(\Omega,\mathbb{R}  ^{d\times d})$, $q$ is the conjugate of $p$, $\mathscr U_\Set$ is the set of all $ u:\Omega\to \bar\Lambda $ measurable  such that  there exists a  $c=c(u,\Omega,\Lambda) >0$  satisfying :
\begin{equation}\label{eq:u_s_intro}
\Big\vert\int_\Omega  u\cdot \text{div}\, \varphi\;dx \Big\vert \leq c \|\varphi\| _{  L^{q}(   {\Omega,   \mathbb{R}^{d\times d}})},\; \quad\forall \varphi\in \Set.
\end{equation}

We point out that  the pseudo-projected gradient  operator depends also on $f$, though the dependence is not exhibited in its notation.
As a first step to approaching \eqref{eq:pro1}, they  have considered  the following perturbed problem: 
\begin{equation} 
\label{eq:pro2}
\inf \left\{ \int_\Omega  \left( f(\nabla_S u)+H(\beta)-F\cdot u \right)dx ;\; (u,\beta)\in \mathscr{U}\right\}.
\end{equation}

The choice of problem \eqref{eq:pro2} is justified by the construction of  a  family of finite dimensional subspaces  $\left\lbrace \Set
_n \right\rbrace_
{n}
$ dense in $W_0^{1,q}(\Omega,\mathbb{R}  ^{d\times d}) $ such that for $u\in \w 1 p{\Omega,\R^d}$, one has
\begin{equation}
\lim_{
{n}
\rightarrow \infty}\intom f(\nabla_{\Set
_n
} u)\;dx=\intom f(\nabla u )\;dx.  \label{eq:approx grad}
\end{equation}
We note that  a $L^p(\Omega,\; \mathbb{R}^{d})-$bounded subset of $\mathscr{U}_\Set$ whose image  by the operator $\nabla_\Set$ is  bounded in $L^p(\Omega,\; \mathbb{R}^{d\times d})$  is not in general strongly pre-compact with respect to the  $L^p(\Omega,\; \mathbb{R}^{d})$ topology.  As a result, compactness of level subsets of the functional in (\ref{eq:pro2}) can not be guaranteed. Nevertheless,    
the authors of  \cite {awi-gan}  have successfully shown existence and, more importantly, uniqueness in \eqref{eq:pro2} %
 under the assumption that 
$F$ is non-degenerate (see definition below). This condition of non-{\tiny }degeneracy  for uniqueness is crucial in a similar problem, the so-called Brenier polar factorization, and more generally, in optimal transport problems. 
Confer \cite{ambrosio2005gradient}, %
\cite{Brenier1991}, %
\cite{gan-ma}, %
\cite{Gangbo1994}, %
\cite{gangbomccann}  %
and
\cite{villani2003topics}.%

In this paper, 
we investigate {the respective roles} played by the strict convexity of $f$, the convexity and smoothness of $H$, and the non-degeneracy of $F$ in problem \eqref{eq:pro2}. More precisely, we impose  less stringent conditions so that  either the map $F$ is allowed to be degenerate or $f$ is allowed to be merely convex or $H$ is neither smooth nor strictly-convex. 
These considerations are not just technicalities. Indeed we note that a prominent case of mere convexity, $ f(\xi)=|\xi| $, is typical for the study of minimal surfaces as well as for the study  of functionals involving the total variation
(see for instance  \cite{Dac_dmcv})
Furthermore, we observe that cases where $H$ is taken to be the characteristic function of a singleton of $\mathbb{R}$ arise in the study of incompressible deformations in Elasticity theory 
(see for instance \cite{gan-wes} and \cite{villani2003topics}).
Finally, the non degeneracy condition tests the extent to which one can hope for uniqueness in the variational problem we considered. 
To deal with these weaker assumptions, we introduce a family of operators $\left\lbrace V_{\Set}^f :  \Set\subset W_0^{1,q}(\Omega,\mathbb{R}  ^d), \;f \;\text{convex} \right\rbrace $  defined by  
\be \label{eq:intom f(nabla_S u) into }
W^{1,p}(\Omega,\mathbb{R}  ^d) \ni u\mapsto V^f_\Set[u]:=\sup_{\phi\in \Set}\intom\left(- u \di \phi-f^*(\phi)\right)\;dx.
\ee
We note that the operator $V_\Set^f$ is actually well defined on the set of measurable functions $u$ defined from $\Omega$ to $ \bar\Lambda$ when the set $\Set$ is a finite dimensional nonempty set and
the function $f$ satisfies appropriate growth conditions.
As a family, these operators  extend the pseudo- projected gradient operators and the distributional gradient. Indeed, $V_\Set^f[u]=\intom f(\nabla_\Set u)$ if $\Set$ is a finite dimensional subspace of $\wo 1 q {\Omega,\R^{d\times d}}$ and $u\in \mathscr U_\Set$  and furthermore $V_\Set^f[u]=\intom f(\nabla u)$ if $\Set=\wo 1 q {\Omega,\R^{d\times d}}$ and $u\in\w 1 p {\Omega,\R^{d}}$. These extensions  are only valid under appropriate conditions on 
 $f$. It is worth pointing out that if $f(\xi)=|\xi|$ and $\Set=\wo 1 q {\Omega,\R^{d\times d}}$ then $V_\Set^f(u)$ is nothing
 but the total variation of $u$ on the set $\Omega$. We show that for a collection of sets $\left\lbrace \Set_n\right\rbrace_{n=1}^{\infty} $ of $W_0^{1,q}(\Omega,\mathbb{R}  ^d)$ satisfying  {\it Hypothesis (H1)} or {\it Hypothesis (H2)} (see section 2),  we have a convergence result in the same spirit as \eqref{eq:approx grad}:
 \begin{equation}
 \lim_{\tau \rightarrow \infty}V^f_{\Set_n}[u]=V^f_{W_0^{1,q}(\Omega,\mathbb{R}  ^d)}[u]\left( = \intom f(\nabla u )\;dx\right)  \label{eq:V_approx grad}
 \end{equation}
 for any $u\in\w 1 p {\Omega,\R^{d\times d}}$ and appropriate conditions on $f$. We thus proceed to study a more general problem :
\begin{equation}
\label{eq:primal H intro}
\inf_{(u,\beta)\in \mathscr U^*_\Set}\left\{ V_\Set^f[u]+\int_\Omega  H(\beta)-F\cdot u\; dx\right\}
\end{equation}
where $\Set$ is an element of a collection of sets satisfying {\it Hypothesis (H1)} or {\it Hypothesis (H2)}, and
\begin{equation}\label{eq:us* intro}
\begin{split}
\mathscr U_\Set^*=\Big\{(u,\beta):& u \in \mathscr U_\Set;\; \beta:\Omega\to[0,\infty);
\\
&\int_\Omega l(u(x))\beta(x)\;dx=\int_\Omega l(y)\;dy \; \forall l \in C_c(\mathbb{R}  ^d)\Big\}.
\end{split}
\end{equation}
Sublevel sets of the integrand in \eqref{eq:primal H intro} are not compact. Nor is  $f$ necessarily strictly convex. However, we show existence and uniqueness in Problem 
\eqref{eq:primal H intro}. In fact, this result holds  for $F$ non-degenerate as well as for a class of degenerate $F$
  provided that the set $\Set$ is chosen 
accordingly ( see  Corollaries \ref{cor:1} and \ref{cor:2}).
Unlike  
optimal transport theory, 
this analysis suggests that  the non-degeneracy condition  is not essential for a uniqueness result in \eqref{eq:pro1}.

Existence and uniqueness results for Problem \eqref{eq:primal H intro} are established thanks to the discovery
of suitable dual problems. Indeed, call $\mathcal C$ the set of all functions $(k,l)$
with
$k,\;l:\mathbb{R}^d  \rightarrow \mathbb{R} \cup\{\infty\}$  Borel measurable, finite at least at one point, and satisfying the relation $l\equiv \infty$ on $\mathbb{R}^d  \setminus \bar \Lambda$ and such that
\begin{equation*}
\label{eq: inequality l, k intro}
k(v) +tl( u) +H(t)\geq u\cdot v\quad \forall u,v\in \mathbb{R}^d  , t>0.
\end{equation*} 
Let $\mathscr A$ be the set of $(k,l,\varphi)$  such that  $(k,l)\in \mathcal{C}$
and $\phi\in  \Set$.
Define the following functional over the set $\mathscr A$:
 $$
 J(k,l,\varphi):=\int_\Omega  k(F+\operatorname{div} \phi)\;dx+\int_\Lambda l\;dy+\intom f^*(\phi)\;dx.
 $$ 
Next, assume that  the map $ F $ and the set $ \Set $ are  such that  for all $ \varphi\in \Set $, 
\begin{equation}\label{eq: Non degenrate cond}
 F+\operatorname{div} \varphi \text{ is non-degenerate.}
\end{equation}
Then $-J$ admits a maximizer $ (k_0,l_0,\varphi_0) $ with $ k_0 $   convex and
{$diam(\Lambda)$-Lipschitz.} As a consequence,
Problem  \eqref{eq:primal H intro} admits a
unique minimizer  $ (u_0,\beta_0) $ and $ u_0 $ satisfies
\begin{equation}\label{eq: EL dual problem}
\begin{cases}
u_0=&\nabla k_0(F+\operatorname{div} \varphi_0)\\
\varphi_0\in&\Phi_\Set(u_0).
\end{cases}
\end{equation}
Here, we have denoted by $\Phi_\Set(u_0)$, the non-empty set of maximizers of problem \eqref{eq:intom f(nabla_S u) into } (see Proposition \ref{pro:convex pseudo projected}).
In order to obtain  condition \eqref{eq: Non degenrate cond}, we consider two distinct situations.

 First, we assume that $ F $ has a countable range, thus degenerate. If  $ \Set $ is an element of a collection of  sets satisfying hypothesis (H2) then it holds that $ F+\operatorname{div} \varphi $ is non degenerate.

Second, we assume $ F $   non-degenerate and  
$ \Set$ is a finite dimensional vector space, as in \cite{awi-gan}. It holds again that $ F+\operatorname{div} \varphi $ is non degenerate.
However, unlike the hypotheses in \cite{awi-gan}, we have allowed the map $f$ to be as singular
as the map $\R ^{d\times d}\ni \xi\mapsto |\xi|$.

We have also studied \eqref{eq:primal H intro}
when $H$ is replaced by 
$H_0:(0,\infty)\to\mathbb{R}  \cup\{\infty\}$ defined by $H_0(1)=0$ and $H_0(t)=\infty$ if $t\neq 1$. 
This case corresponds to the case of measure preserving maps.
Note that $H_0$ is not even continuous. However, it may be obtained as a limit of functions $H_n$ which are $C^1(0,\infty)$ convex functions and satisfy \eqref{eq:cond}. 
We show that for such singular $H_0$, the  corresponding problem
\begin{equation}
\label{eq:primal H0 intro}
\inf_{u\in \mathscr U^1_\Set }\left\{ V_\Set^f[u]-\int_\Omega  F\cdot u \;dx\right\}
\end{equation}
with
\begin{equation}\label{eq:us1 intro }
\mathscr U_\Set^1=\left\{u\in \mathscr U_\Set: \int_\Omega l(u(x))\;dx=\int_\Omega l(y)\;dy \; \forall l \in C_c(\mathbb{R}  ^d)\right\}
\end{equation}
admits a unique minimizer. (See Theorem \ref{thm:main 2}).

To obtain existence and uniqueness results in
problem \eqref{eq:primal H0 intro}, we exploit a dual formulation 
and maximize $-J$ over the set that consists of 
$(k,l,\varphi)$   such that
$\phi\in \Set$ and   $k,\;l:\mathbb{R}^d  \rightarrow \mathbb{R} \cup\{\infty\}$ are  Borel measurable, finite at least at one point, and satisfy the relations $l\equiv \infty$ on $\mathbb{R}^d  \setminus \bar \Lambda$ and 
\begin{equation*}
\label{eq: inequality l, k intro 1}
k(v) +l( u) \geq u\cdot v\quad \forall u,v\in \mathbb{R}^d.
\end{equation*} 
One shows that $-J$
admits a maximizer $ (k_0,l_0,\varphi_0) $ with $ k_0  $ convex and Lipschitz and the unique minimizer of
problem  \eqref{eq:primal H0 intro} is  $ u_0 $ given by 
$$
u_0=\nabla k_0(F+\operatorname{div} \varphi_0).
$$
Finally, we show  convergence of a sequence of  problems of the form \eqref{eq:primal H intro} to \eqref{eq:primal H0 intro}.
More precisely, we show 
that the minimizer of problem \eqref{eq:primal H0 intro} may be obtained as
limit of minimizers of problems of the form 
\eqref{eq:primal H intro} provided that the dual problems admit  regular enough maximizers.
In fact, suppose the map $ F $ and the set $\Set$ are  such that  for all $ \varphi\in \Set $,  the map $ F+\operatorname{div} \varphi $ is non-degenerate.
For $(u,\beta) \in \mathscr U_\Set$, define
\[
 I_n(u,\beta)=V_\Set^f[u]+\intom H_n(\beta)-u\cdot F\;dx
 \label{eq:primal n intro}
\]
and set 
\[
 I_0(u )=V_\Set^f[u]-\intom u\cdot F\;dx.
 \label{eq:primal 0 intro}
\]
Thanks to Theorem \ref{thm:main 1}, the problem
\be\label{eq:approx n}
 \inf_{(u,\beta)\in \mathscr U_\Set^*  }I_n(u,\beta)
\ee
admits a unique minimizer that we denote 
$(u_n,\beta_n)$ with $u_n=\nabla k_n(F+\di \phi_n)$ for some $k_n:\R^d\to \R$ 
convex and $\phi_n\in \Set$.
Denote $u_0$ the unique minimizer of \eqref{eq:primal H0 intro}. 
If for all $n\in \N^*$ the map $k_n$ is differentiable then 
the sequence $\{ u_n \}_{n\in \N^*}$ 
 converges almost everywhere to   $u_0$  and in addition, the minima $\{  I_n(u_n,\beta_n)\}_{n\in \N^*}$ 
   converge to $I_0(u_0)$  (Cf. Theorem \ref{thm:main 3}).

\section{Preliminaries}
\subsection{\bf Notation and definitions.}
\begin{itemize}%
\item[$ \bullet $] Throughout this manuscript, $ \Omega$ and $\Lambda\subset\mathbb{R}^d$  are two bounded open convex sets;
$r^*>1$ is  such 
that $B(0,1/r^*)\p \Lambda\subset   B(0,r^*/2)$;  $p\in (1,\infty)$ and $q$ is its conjugate, that is, 
 $p^{-1}+q^{-1}=1$.
\item[$ \bullet $] Given $A\subset \R^d$, the indicator function  of $A$ is defined as $$\chi_{A}(x)=\begin{cases}0 \qquad \text{ if }  x\in A, \\ \infty\qquad \text{ otherwise.}\end{cases}
 $$
\item[$ \bullet $] For any  subset $\mathcal{S}$ of $ W_0^{1, q}( {\Omega,\mathbb{R}   ^{d\times d}})$, we denote by $\text{span}(\Set)$ the linear subspace  of 
$W_0^{1, q}( {\Omega,\mathbb{R}   ^{d\times d}})$   generated by $\Set$.

 \item[$ \bullet $] We denote by $f^*$ the Legendre transform of  the map $f:\mathbb{R}   ^{d\times d}\longrightarrow \mathbb{R}   $ so that %
 \[
 f^*(\xi^*)=\sup_{\xi\in \mathbb{R}   ^{d\times d}}\left\{\xi\cdot \xi^*-f(\xi)\right\}.
 \] 
 \item[$ \bullet $] If $ h : \R^d\longrightarrow \R\cup \{\infty\}$ is convex then the subdifferential $\partial h (x)$ of $h$ at $x\in \text{Dom} (h)$ is closed and convex. If $\partial h (x)$  is  non-empty
 we denote by $grad[h](x)$ the element of $\partial h (x)$ with minimum norm :
$$
|grad[h](x)|= \min \left\lbrace |y|: y\in\partial h (x)\right\rbrace; \qquad x\in \text{Dom}(h).
$$
 
\item[$ \bullet $]  Let $\Set\subset W_0^{1, q}( {\Omega,\mathbb{R}   ^{d\times d}})$. We
  denote by  $\mathscr S_f$ the set 
 \begin{equation}
 \mathscr S_f :=\left\{\phi\in \Set: \int_\Omega f^*(\phi) \text{ is finite}\right\}.
 \end{equation}
 \item[$ \bullet $] Let $F :\mathbb{R}^d \longrightarrow \mathbb{R}^d $  be measurable. We say that $F$ is non-degenerate if  for any $N\p\mathbb{R}^d  $ such that $\mathcal{L}^d(N) =0 $ we have $\mathcal{L}^d(F^{-1}(N)) =0. $
 
\end{itemize}

\subsection{\bf Assumptions.}

\begin{itemize}%
\item[\bfseries{(A0)}] We additionally assume  that 
there exists a strictly convex function that is $ C^1(\bar \Omega) $ and vanishes on the boundary of $ \Omega $.
\item[\bfseries{(A1)}] The set 
$\Set$ is a subset of $W_0^{1, q}( {\Omega,\mathbb{R}   ^{d\times d}})$. In addition, the map $f:\R^{d\times d}\to \R$ is convex and satisfies the following three properties:
\begin{itemize}
\item[(i)] There exist  $a,b,c >0$ such that for all $\xi\in \R^{d\times d}$,  
    \begin{equation}\label{eq:f growth}
    c\frac{|\xi|^p}p+b\geq f(\xi)\geq a|\xi|-b
    \end{equation}
and
for all $\xi^*\in \pd f(\xi) $,
\begin{equation}\label{eq:pd f growth}
    |\xi^*|^q\leq  c|\xi|^p + b.
\end{equation}
 \item[(ii)]  
 The set $\mathscr S_f$ is non empty. 
\item [(iii)] One of the following two conditions holds:
\begin{enumerate}
 \item[(a)] The map $f$ is such that  $\partial f^*(x^*)$ is non-empty and $grad[f^*](x^*)=0$
for each $x^*\in \text{Dom} f^*$.
\item[(b)] The map $f$ is 
strictly convex and there exist $\bar a,\bar b>0$ such that for all $\xi^*\in \mathbb{R}^{d\times d}$,
one has
 \begin{equation}\label{eq: growth on f* and grad f*}
 f^*(\xi^*)\leq \bar a+\bar b|\xi^*|^{q}
\quad\hbox{and}\quad  |\nabla f^*(\xi^*)|\leq \bar a +\bar b|\xi^*|^{q-1}.
 \end{equation}
 
\end{enumerate}   
  
\end{itemize}

 \item[\bfseries{(A2)}] The map $H$ is $C^1(0,\infty)$, strictly convex, and  such that  
 \[
 \lim_{t\rightarrow 0^+}   H(t)= \lim _{t\rightarrow \infty} \frac{H(t)}{t}=+\infty.
 \]
 
 \item[\bfseries{(A3)}] The function $F$ is measurable and belongs to $ L^1(\Omega)$.
 \\
 Let $ \Set$ be a subset of $ W_0^{1,q }({\Omega,\mathbb{R}^{d\times d}})$. We say that $F$ satisfies the condition $\textbf{(ND)}_{\Set}$ if
 $$   \di(\phi) + F \;\text{ is non degenerate}     $$
 for all $\phi\in \Set$.
\end{itemize}

\begin{rmk} \label{rem:pteOf_h}

\begin{itemize}
\item[(i)]  As  $f$ satisfies (\ref{eq:f growth}), we have
\be \label{eq:f^* growth}
-b+c^p\frac{|\xi^*|^q}{q}\leq f^*(\xi^*)\leq \chi_{\bar B(0,a)}(\xi^*)+b
\ee
for all $\xi^*\in \R^{d\times d}$. 
\item[(ii)] If $f$ satisfies case (b) in (iii)  of Assumption $\mathbf{(A1)}$,
then $f^*$ is continuously differentiable. In that case, $grad[f^*]=\nabla f^*$. 
\item[(iii)] If $f$ satisfies case (a)  of Assumption $\mathbf{(A1)}$(iii) then $ 0\in \partial f^*(x^*) $ for every element $ x^* \in \operatorname{Dom}(f^* )$. Consequently, the map $ f^* $ is constant on $\operatorname{Dom}(f^* )$ and the following equalities are satisfied for all $ x^* $ and $ y^* $ in $\operatorname{Dom}(f^* ) $:
\begin{equation}\label{eq:f* const}
f^*(x^*)-f^*(y^*)=grad[f^*](x^*)=grad[f^*](y^*)=0.
\end{equation}
 
\item[(iv)] The assumption $\mathbf{(A0)}$  is satisfied by $\Omega=B(0,1)\subset \mathbb{R}^d$ with the strictly convex function being the map $\mathbb{R}^d\ni x\mapsto |x|^2-1 $.
\item[(v)] The map $f:\R^{d\times d}\to \R$ 
defined by $f(\xi)=|\xi|$ 
satisfies case (a) in (iii)  of Assumption $\mathbf{(A1)}$.
The map $f:\R^{d\times d}\to \R$ 
defined by $f(\xi)=|\xi|^p$ 
satisfies case (b) in (iii)  of Assumption $\mathbf{(A1)}$.

\end{itemize}
\end{rmk}

The following Lemma summarizes some elementary properties of  $ H $. We refer the reader to Remark 2.1 in \cite{awi-gan}. 
\begin{lemme} \label{lem:pteOf_h}
 Assume 
$\mathbf{(A2)}$ 
holds. Then,
 \begin{enumerate}
 \item[(i)] The map $H':(0,\infty)\to \mathbb{R}  $ is a strictly increasing  bijection.
  \item [(ii)]The Legendre transform $H^*$ of $H$ is a strictly increasing bijection from $\mathbb{R}  $ to $\mathbb{R}  $.\label{lem_part:h*inc}  
  \item[(iii)] \label{lem_part:h*_growth}
Let  $g:\mathbb{R}  \to \bar{\mathbb{R} } $
be defined by $g(s)=\alpha s-\beta H^*(s)$, with $\alpha,\beta>0$.
Then
$$\lim_{s\to -\infty}g(s)=\lim_{s\to \infty}g(s)=-\infty.$$
\end{enumerate}
\end{lemme}
Define $H_0$ by 
 \begin{equation} \label{eq H_0}
 H_0(t)=\begin{cases}
 0\qquad \;\; t=1\\
 \infty\qquad t\neq 1\\
 \end{cases}
 \end{equation}
 and, for $n\geq 1,$ 
 \begin{equation} \label{eq H_n}
 H_n(t)= H(t)- H(1)+ n(t-1)^2. 
 \end{equation}
  The following Lemma is straightforward.
 \begin{lemme} \label{lem: lim H0}
   Assume  
$\mathbf{(A2)}$ 
holds. Then,
  \begin{itemize}
   \item[(i)]  There exists $\bar H\in \R$ such that
   $$ \bar H= \min_{t\in[0,\infty)} H(t).$$ 
   \item [(ii)] The collection  $\left\lbrace H_n\right\rbrace_{n=1}^{\infty} $ is a non decreasing  sequence of functions that converges pointwise to $H_0$. In addition, for all $n\in \N^*$, %
   the map $H_n$ is a $C^1(0,\infty)$ strictly convex  function that satisfies 
   \begin{equation*} 
   \lim _{t\rightarrow 0^+}  H_n(t)= \lim _{t\rightarrow \infty} \frac{H_n(t)}{t}=+\infty.
   \end{equation*}
   \item [(iii)] Let $t>0$. If  $\left\lbrace H_n(t)\right\rbrace_{n=1}^{\infty} $ is uniformly bounded above by a constant $c_0$ then
   $$ n(t-1)^2 \leq  c_0 +H(1)-\bar H $$
   and $t=1$.
  \end{itemize}

 \end{lemme}

\subsection{\bf Hypothesis on the underlying sets of pseudo-gradients. }

 We recall that in \cite{awi-gan}, %
the construction of $\nabla_{\Set^{\tau}} u$ has relied on hypothesis on
the underlying sets $ \Set^{\tau} $ that we summarize in {\it Hypothesis (H1)} below. 
\\ 
{\it Hypothesis (H1).}

A collection
$  \{ \mathfrak{A}_{n }\}^{\infty }_{ n=1} $ of subsets of $ W_0^{1,q }({\Omega,\mathbb{R}^{d\times d}}) $  satisfies   {\it Hypothesis (H1)}  if

\begin{itemize}
 \item[(i)] $ \mathfrak{A}_{n }$ of a finite dimensional subspace of  $ W_0^{1,q }({\Omega,\mathbb{R}^{d\times d}}) $ for each $n\in \mathbb{N^{*}}$.
 \item[(ii)]  The map $ \nabla \varphi $ has a countable range whenever
 $ \varphi\in \mathfrak A_n $, for any $ n\in \mathbb{N}  ^* $.
 \item[(iii)] The set $ \cup _{n\in \mathbb{N}  ^*}\mathfrak{A} _{n} $ is dense in 
 $ W_0^{1,q }({\Omega,\mathbb{R}   ^{d\times d}} )$.
 
 \item[(iv)] For  $i\leq j$, we have the inclusion 
 $\mathfrak{A}_{i }\subset \mathfrak{A}_{j}.$
\end{itemize}
An explicit construction of sets satisfying {\it Hypothesis (H1)}   is provided in \cite{awi-gan}. Here, we build on the conditions of {\it Hypothesis (H1)}   and  we relax conditions on the underlying sets:  
\\
\noindent
{\it Hypothesis (H2).} 

A collection
$  \{ \mathfrak{Q}_{n }\}^{\infty }_{ n=1} $ of subsets of $ W_0^{1,q }({\Omega,\mathbb{R}^{d\times d}}) $  satisfies  {\it Hypothesis (H2)} if
\begin{itemize}
 \item[(i)]  $\text{Span}(\mathfrak{Q}_{n})$ is of finite dimension and $ \mathfrak{Q}_{n }$ is a non-empty closed and convex  subset of $ W_0^{1,q }({\Omega,\mathbb{R}^{d\times d}})$.
\item[(ii)] The map $ \operatorname{div} \varphi $ is non-degenerate whenever
$ \varphi\in \mathfrak{Q}_{n} $, for any  $ n\in \mathbb{N}  ^* $.
 \item[(iii)] The set $ \cup _{n\in \mathbb{N}  ^*}\mathfrak{Q} _{n} $ is dense in 
 $ W_0^{1,q }({\Omega,\mathbb{R}   ^{d\times d}} )$.%
  
 \item[(iv)]  For  $i\leq j$, the inclusion
 $%
 \mathfrak{Q}_{i }\subset \mathfrak{Q}_{j}
 $ %
 holds.

\end{itemize}
The next lemma asserts that a collection of sets can be constructed to satisfy {\it Hypothesis (H2)}. 
\begin{lemme} \label{ lem : existence H2}
Assume $\mathbf{(A0)}$ holds. Then, there exists a collection of sets\newline $  \{\mathfrak{Q}_{n }\}^{\infty } _{n=1}$ satisfying the requirements of {\it Hypothesis (H2)}. 
\end{lemme}
\begin{rmk}
The condition $\mathbf{(A0)}$ in Lemma \ref{ lem : existence H2} is only needed for requirement (ii) of {\it Hypothesis (H2)}. 
\end{rmk}
{Proof.}
Suppose that 
$ \psi $ is a strictly convex function that is $ C^1(\bar \Omega) $ and vanishes on the boundary of $ \Omega $ as given by Assumption \textbf{(A0)}.
Let $ \phi_0:\Omega\to \R^{d\times d} $ be defined by
\[\phi_0=
\begin{pmatrix}
\psi&0&\cdots&0
\\
0&\psi&\cdots&0
\\
\vdots&\vdots&\ddots&\vdots
\\
0&0&\cdots&\psi
\end{pmatrix}.
\]
As $\psi$ is $C^1(\bar \Omega) $, we have $\phi_0\in W_0^{1,q }({\Omega,\mathbb{R}^{d\times d}}) $
and it follows that $ \di \phi_0=\nabla \psi $. Thus, for almost every $ x $ in $ \Omega $, we have
\[
\det(\nabla( \di \phi_0)(x))=\det(\nabla^2 \psi(x))>0.
\]
Thanks to Lemma 5.5.3 in \cite{ambrosio2005gradient}, the map
$ \di \phi_0 $ is non-degenerate.
Let $  \{ \mathfrak{A}_{n }\}^{\infty }_{ n=1} $ be a collection of sets satisfying Hypothesis (H1).
One readily checks that 
the family of sets defined by
\[
\mathfrak{Q} _{n}= \left \{ \varphi+\epsilon \varphi_0: \varphi\in \mathfrak{A}_n;  \epsilon\geq \frac1n \right \}
\]
for $n\in \mathbb{N}  ^*$, satisfies hypothesis (H2).\endproof \cqfd
\subsection{\bf Special displacements.}
 To $\Set\subset\wo{1}{q}{\Omega,\mat dd}$ we associate $ \mathscr U_\Set $, the set of all 
 $ u:\Omega\to \bar\Lambda $ measurable  such that  there exists  $\bar c=\bar c(u,\Omega,\Lambda) >0 $  satisfying :
\begin{equation}\label{eq:u_s}
\Big\vert\int_\Omega  u\cdot \text{div}\, \varphi \;dx\Big\vert \leq\bar  c \|\varphi\| _{  L^{q}(   {\Omega,   \mathbb{R}^{d\times d}})}\; \quad\forall \varphi\in \Set.
\end{equation}
Remark that if $u \in \mathscr U_\Set$, then $ u $ belongs to $ L^\infty(\Omega,\mathbb{R}^d) $ since $u$ has values in $ \bar \Lambda $ which is bounded.  
If $ \text{span}(\Set)$ is of finite dimension then  $\mathscr U_\Set  $ is the set of all measurable maps $u:\Omega\to \bar\Lambda$. 
In fact, the linear map $ \text{span}(\Set)\ni \varphi\mapsto \int_\Omega u\operatorname{div}\varphi $ is continuous with respect to the $ L^q $- norm as in finite dimension, all norms are equivalent. Therefore, we may find $c$
for which Equation \eqref{eq:u_s} holds for all 
$\varphi\in \text{span}(\Set)$ and in particular
for all 
$\varphi\in \Set$.
\\
At any rate, $\mathscr U_\Set$ contains $\w{1}{p}{\Omega,\R^d}$. Indeed, notice that for  a fixed
$u\in \w{1}{p}{\Omega,\R^d}$, we have, for all 
$\phi\in \Set$:
\[\left|\int_\Omega u\cdot \operatorname{div}\varphi\;
dx\right|=\left|-\int_\Omega \langle\nabla u,\varphi\rangle \;dx\right|
\leq  \left\|\nabla u\right\|_{ L^{p}(   {\Omega,   \mathbb{R}^{d\times d}})}
\left\| \varphi\right\|_{ L^{q}(   {\Omega,   \mathbb{R}^{d\times d}})}.
\]
We introduce the following set
\begin{equation*}%
\mathscr U_\Set^1=\left\{u\in \mathscr U_\Set: \int_\Omega l(u(x))\;dx=\int_\Lambda l(y)\;dy \; \forall l \in C_c(\mathbb{R}  ^d)\right\}
\end{equation*}
and
\begin{align*}%
\mathscr U_\Set^*=\Bigg\{(u,\beta):& u \in \mathscr U_\Set;\; \beta:\Omega\to[0,\infty);\\ &\int_\Omega l(u(x))\beta(x)\;dx=\int_\Lambda l(y)\;dy \; \forall l \in C_c(\mathbb{R}  ^d)\Bigg\}.
\end{align*}

Notice  that $\mathscr U_\Set^1=\{u\in \mathscr U_\Set: (u,1)\in \mathscr U_\Set ^*\}$. This corresponds to measure preserving displacements.

\subsection{\bf Extended pseudo-projected gradient}
Let   $\Set\subset\wo{1}{q}{\Omega,\mat dd}$  and  $ u \in \mathscr U_\Set$. Define
 \[
 \mathcal G_\Set( u ):=
 \left\{
 G\in \el p {\Omega,\R^{d\times d}}:
 \intom  u \di \phi\;dx =-\intom \pair{G}{\phi}\;dx\;
 \forall \phi\in \Set  
 \right\}.
 \]
Consider the operator
\be \label{eq:intom f(nabla_S u) }
V^f_\Set
( u ):=\sup_{\phi\in \Set}\intom\left(- u \di \phi-f^*(\phi)  \right)\;dx =\sup_{\phi\in \mathscr S_f}\intom\left(- u \di \phi-f^*(\phi)  \right)\;dx.
\ee
We denote by $\Phi_S( u )$ the   set of maximizers of Problem \eqref{eq:intom f(nabla_S u) }.
\begin{lemme}
Let $\Set\subset\wo{1}{q}{\Omega,\mat dd}$  and  $ u \in \mathscr U_\Set$.
\begin{enumerate}
 \item We have
 \[
 \mathcal G_\Set( u )=
 \left\{
 G\in \el p {\Omega,\R^{d\times d}}:
 \intom  u \di \phi dx=-\intom \pair{G}{\phi}dx;\ 
 \forall \phi\in \operatorname{span}(\Set)  
 \right\}.
 \]
 \item If $\operatorname{span}(\Set)$ 
 is finite dimensional, then 
 $\mathcal G_\Set( u )$ is
 nonempty.
 
\end{enumerate}

\end{lemme}
Proof. Set 
\begin{align*}
 \mathcal{\bar G}_\mathcal{S}(u)&=\left\{
 G\in L^p({\Omega,\mathbb{R}^{d\times d}}):
 \int_\Omega   u \operatorname{div}\phi dx=-\int_\Omega \langle G,\phi\rangle dx \;
 \forall \phi\in \operatorname{span}(\mathcal S) 
 \right\}.
\end{align*}
As $\mathcal{S}\subset \operatorname{span}(\mathcal{S})$,
we have  
$\mathcal{\bar G}_\mathcal{S}(u)
\subset 
 \mathcal{ G}_\mathcal{S}(u)
 $.
Next, let $G\in  \mathcal{ G}_\mathcal{S}(u)$. Assume that 
$\phi\in \operatorname{span}(\mathcal{S})$.
We may find $n\in \mathbb N$, $\lambda_1,\dots,\lambda_n\in \R$
and $\phi_1,\dots,\phi_n\in \mathcal{S}$
such that
$\phi=\sum_{i=1}^n\lambda_i\phi_i.$
Then
\[
\int_\Omega u\operatorname{div}\phi dx =
\int_\Omega u\operatorname{div}\sum_{i=1}^n\lambda_i\phi_i dx
=\sum_{i=1}^n \lambda_i\int_\Omega u\operatorname{div}\phi_i dx
=\sum_{i=1}^n -\lambda_i\int_\Omega \langle G,\phi_i\rangle dx
\] 
and
\[-\int_\Omega \langle G,  \phi\rangle dx
=-\int_\Omega \langle G,\sum_{i=1}^n \lambda_i\phi_i\rangle dx
=\sum_{i=1}^n -\lambda_i\int_\Omega \langle G,\phi_i\rangle dx.
\] 
Thus $G\in\mathcal{\bar  G}_\mathcal{S}(u)
 $. We deduce that $\mathcal{G}_\mathcal{S}(u)
\subset 
 \mathcal{\bar  G}_\mathcal{S}(u)
 $.
It follows that part 1.) holds. To obtain part 2.), we use part 1.) and the Riesz Representation Theorem.
\cqfd

The following results are  essentially found in  Proposition 3.1 in \cite{awi-gan}.%

\begin{pro} \label{lem:pseudo from awi-gan}
Suppose   that the set  $\Set$ is a  finite dimensional subspace of \\
$ W_0^{1,q }({\Omega,\mathbb{R}^{d\times d}})$  and $f$ is  $C^1$ 
 and strictly convex. Suppose, in addition that there exist constants $c_1,c_2,c_3>0$ such that
\begin{align*}%
 -c_3+c_2|\xi|^p\leq f(\xi)\leq c_3+c_1|\xi|^p\\
   |Df(\xi)|\leq c_3+c_1|\xi|^{p-1}\\
 |Df^*(\xi)|\leq c_3+c_1|\xi|^{q-1}
\end{align*} 
for all $\xi\in \mathbb{R}  ^{d\times d}$. Then, there exists a unique map denoted $\nabla_\Set  u $ that minimizes
 $$\inf_{G\in \mathcal G_\Set( u )}\intom f(G)\;dx.$$ 
 Moreover, $\nabla_\Set u $ uniquely satisfies
 $G\in \mathcal G_\Set( u )$ and  $ Df(G)\in \Set. $
\end{pro}

\noindent In the next Proposition, we establish  similar results as in Proposition \ref{lem:pseudo from awi-gan} but under weaker  assumptions on $\Set$ and $f$ (except in part 4).
\begin{pro}
 \label{pro:convex pseudo projected} Assume  $\mathbf{(A1)}$  holds. Assume  $\Set$ is a finite dimensional non-empty closed and convex  subset of $ W_0^{1,q }({\Omega,\mathbb{R}^{d\times d}})$  and let $ u \in \mathscr U_\Set$.
\begin{enumerate} 
 
\item For all $G\in \mathcal G_\Set( u )$, $\phi\in \Set$, we have
\[ %
\intom f(G) \;dx\geq \intom\left(- u \di \phi-f^*(\phi)  \right)\;dx.
\]

\item The supremum in problem 
\eqref{eq:intom f(nabla_S u) }
is attained. 

 \item
 A map $\bar \phi$ belongs to $\Phi_\Set( u )  $ if and only if $\bar\phi$ belongs to $\mathscr S_f $ 
 and
 \[
  \intom \left( grad[f^*](\bar\varphi)\cdot \left(\phi-\bar \phi\right)
  + u  \cdot \left(\di\phi- \di\bar\phi\right) \right)\;dx\geq 0
 \]
  for all $\phi\in\mathscr S_f $.

 \item Suppose that  the hypotheses of Proposition \ref{lem:pseudo from awi-gan} are satisfied.
 Then we have 
 $$\intom f(\nabla_\Set  u )\;dx=V_\Set^f( u )$$ 
 and $\Phi_\Set( u )=\{Df(\nabla_\Set  u )\}$.%
\end{enumerate}

\end{pro}

\noindent 
Proof. 
1.) Let $\varphi\in \Set$ and $G\in \mathcal G_\Set( u )$, By using the Legendre transformation,
$$
\intom f(G)\;dx\geq \intom G\cdot \phi -f^*(\phi)\;dx= \intom - u \cdot \di\phi -f^*(\phi)\;dx. 
$$
2.) Let  $\phi\in \Set$. We use (\ref{eq:u_s}) and (\ref{eq:f^* growth}) to get  
\begin{equation}
\label{eq: coercivity T}
\begin{aligned}
\intom\left( u \di \phi+f^*(\phi)  \right)\;dx\geq &  -\bar c \|\phi\|_{\el q {\Omega,\R^{d\times d}}}+ \intom f^*(\phi )\;dx\\
\geq& -\bar c \|\phi\|_{\el q {\Omega,\R^{d\times d}}}
+q\s c^{-q} \|\phi\|^q_{\el q {\Omega,\R^{d\times d}}}.
\end{aligned}
\end{equation}
In light of (\ref{eq: coercivity T}), $q>1$ implies that
 the map $$\mathscr S_f\ni  \phi\mapsto T(\phi):=\intom\left( u \di \phi+f^*(\phi)  \right)\;dx$$  is $L^q-$ coercive. 
Moreover, the convexity of $f^*$ guarantees that $T$  is  lower semi-continuous.
 The direct methods of the calculus of variations thus yield the existence of a
maximizer in problem \eqref{eq:intom f(nabla_S u) }.
\\
3.)
Let $\bar \phi\in\Phi_\Set( u )  $  so that  $\bar\phi\in\mathscr S_f$. Let  $\phi\in\mathscr S_f$ and $\epsilon\in (0,1)$. The convexity of $f^*$  ensures that  $\bar \phi+\epsilon  (\phi-\bar \phi)\in\mathscr S_f$ and the maximality property of  $\bar \phi $ implies that 
\begin{equation}
\label{eq: maximalty property of phi 1}
  \intom u \cdot \di \bar \phi +f^*(\bar \phi)\;dx\leq 
  \intom  u \cdot \left(\di \bar \phi+\epsilon\di (\phi-\bar \phi)  \right)\;dx +
  f^*(\bar \phi+\epsilon  (\phi-\bar \phi))\;dx.
   \end{equation}
We rewrite (\ref{eq: maximalty property of phi 1}), in turn, as
\begin{equation}
\label{eq: maximalty property of phi 2}
\intom \frac{ f^*(\bar \phi+\epsilon  (\phi-\bar \phi))-f^*(\bar \phi)}{\epsilon}+  u \cdot \di (\phi-\bar \phi)\; dx
\geq 0.
\end{equation}
Note that $ grad[f^*](\bar\varphi+\epsilon(\varphi-\bar\varphi)) $ belongs to the set 
$ \partial f^*((\bar\varphi+\epsilon(\varphi-\bar\varphi))) $ whenever $(\bar\varphi+\epsilon(\varphi-\bar\varphi)) $ is in the domain of $ f^* $. It follows that

\[\int_\Omega\left(grad[f^*]((\bar\varphi+\epsilon(\varphi-\bar\varphi)))\cdot(-\epsilon(\varphi-\bar\varphi))\right)\;dx\leq \int_\Omega\left(f^*(\bar\varphi)-f^*(\bar\varphi+\epsilon(\varphi-\bar\varphi))\right)\;dx\]
that is,
\begin{equation}\label{eq: EL 0}
\int_\Omega\left(grad[f^*](\bar\varphi+\epsilon(\varphi-\bar\varphi))\cdot(\varphi-\bar\varphi)\right)dx\geq 
\int_\Omega\frac{f^*(\bar\varphi+\epsilon(\varphi-\bar\varphi))-f^*(\bar\varphi)}
{\epsilon}\;dx.
\end{equation}
 We combine \eqref{eq: maximalty property of phi 2} and \eqref{eq: EL 0} to get
\begin{equation}\label{eq:gradf^*beforelimit} 
\int_\Omega\left(grad[f^*](\bar\varphi+\epsilon(\varphi-\bar\varphi))\cdot( \varphi-\bar\varphi)+u\operatorname{div}(\varphi-\bar\varphi)\right)dx\geq 0.
\end{equation}
First, we assume that $\mathbf{(A1)}$(iii)(a) holds. In light of \eqref{eq:f* const}, we have $grad[f^*](\bar\varphi+\epsilon(\varphi-\bar\varphi)=grad[f^*](\bar\varphi)$.  Equation
\eqref{eq:gradf^*beforelimit}
becomes 
\[
\int_\Omega\left(grad[f^*](\bar\varphi )\cdot( \varphi-\bar\varphi)+u\operatorname{div}(\varphi-\bar\varphi)\right)dx\geq 0 .\]
Second, we assume that $\mathbf{(A1)}$(iii)(b) holds.
In light of {Remark \ref{rem:pteOf_h} (ii)}, 
we use the growth condition on $ \nabla f^* $ in \eqref{eq: growth on f* and grad f*}, the Lebesgue dominated convergence theorem  and let $\epsilon$ go to $0$ in \eqref{eq:gradf^*beforelimit} to  obtain that:
\[ 
\int_\Omega\left(grad[f^*](\bar\varphi )\cdot( \varphi-\bar\varphi)+u\operatorname{div}(\varphi-\bar\varphi)\right)dx\geq 0.\]

We next show the converse implication. Let $\varphi\in \mathscr S_f $ such that
\begin{equation}\label{eq:EL1}
0\leq 
\int_\Omega\left(u\operatorname{div}(\varphi-\bar\varphi)+grad[f^*](\bar\varphi )\cdot( \varphi-\bar\varphi)\right)dx,
\end{equation}
for all $\varphi\in \mathscr S_f$. We notice that, as $f^*$ is convex, the range
 of the map $grad[f^*](\bar \phi)$ lies in the sub-differential of $f^*$ so that $ f^*(\phi)-f^*(\bar\phi)\geq grad[f^*](\bar \phi)\left(\phi-\bar \phi \right)$ for all $\varphi\in \mathscr S_f$. Then, the inequality \eqref{eq:EL1} implies that
\begin{align*}
0\leq & \int_\Omega\left(u\operatorname{div}(\varphi-\bar\varphi)+ (f^*(\varphi)-f^*(\bar\varphi ))\right)dx
\end{align*}
for all $\varphi\in \mathscr S_f$, that is, 
\begin{align*}
\int_\Omega\left(u\operatorname{div}\bar\varphi+  f^*(\bar\varphi) \right)dx\leq & \int_\Omega\left(u\operatorname{div} \varphi + f^*(\varphi )\right)dx
\end{align*} 
for all $\varphi\in \mathscr S_f$. We conclude that $\bar\varphi\in \Phi_\Set( u ) $.\\
4.) Thanks to Proposition \ref{lem:pseudo from awi-gan}, $D f(\nabla_\Set u)\in \Set$. Next, we set $\phi_0:= D f(\nabla_\Set u)$. By definition of $f^*$,
$$
 f(\nabla_\Set  u ) + f^*(\phi )\geq  \phi \cdot \nabla_\Set  u
$$
for all $\phi\in \Set$. As $f$ is convex and $\phi_0 = D f(\nabla_\Set u)$,  
we have
$$
f(\nabla_\Set  u )+ f^*(\phi_0)  =\phi_0\cdot \nabla_\Set  u,
$$
Thus,
$$
\intom f(\nabla_\Set  u )\;dx\geq \intom \phi \cdot \nabla_\Set  u \;dx-\intom f^*(\phi )\;dx=
\intom - u \di\phi\;dx -\intom f^*(\phi )\;dx
$$
and
\[
 \intom f(\nabla_\Set  u )\;dx=\intom \phi_0\cdot \nabla_\Set  u\;dx -\intom f^*(\phi_0)\;dx=
\intom - u \di\phi_0 \;dx-\intom f^*(\phi_0)\;dx.%
\]

We deduce that $\phi_0\in \Phi_\Set( u )$.
Since   $f^*$ is strictly convex, we conclude that $ \Phi_\Set( u )=\{D f(\nabla_\Set u)\}$
and moreover, $\intom f(\nabla_\Set u )=V^f_\Set( u )$.
\cqfd 
In the next Proposition, we establish a convergence result in the spirit of \eqref{eq:approx grad}. We also connect the operator $V_\Set^f$ with the usual notions of gradient and total variation.

\bp \label{lem:vsf and convergence}
Assume  $\mathbf{(A1)}$  holds.  Assume that $ \Set_n$ is a finite dimensional non-empty closed and convex  subset of $ W_0^{1,q }({\Omega,\mathbb{R}^{d\times d}})$  for each $n\geq 1$. The following holds.
\begin{enumerate}
\item If $ \left\lbrace \Set_n\right\rbrace_{n=1}^{\infty} $ is a monotonically increasing family of subsets of some set $ \Set_0 $ and $ \cup _{n\in \N^*}\Set_n $ is dense in $ \Set_0 $ with respect to the $ \wo 1 q {\Omega,\R ^{d\times d}} $ norm then
\[
\Lii n V_{\Set_n}^f[u]=V_{\Set_0}^f[u]
\]
for any  
$ u\in\mathscr U_{\Set_0}. $
\item If $ \Set=\wo 1 q {\Omega,\R^{d\times d} }$  and $ u\in \w 1 p {\Omega,\R^{d} }$ then $ V_\Set^f[u]=\intom f(\nabla u)dx $.

\item Assume $ u\in BV(\Omega,\R^{d\times d}) $ and $ f(\xi)=|\xi| $ for all $ \xi\in  \R^{d\times d}$. If\\
$ \Set=\wo 1 q {\Omega,\R^{d\times d}} $
then
$ V_\Set^f[u] $ is the total variation of $ u $.
\end{enumerate}
\ep

\begin{rmk}
A consequence of   Proposition \ref{lem:vsf and convergence} is the following : If the sequence of sets $ \{\Set_n \}_{n\in \N^*}$ is monotonically increasing to $ \wo 1 q {\Omega,\R ^{d\times d}}$ and $ u\in \w 1 p {\Omega,\R^{d} }$ we have
\[
\Lii n V_{\Set_n}^f[u]=\intom f(\nabla u)\;dx.
\]
\end{rmk}

\proof 
1.) Recall that
\[
V_{\Set_n}^f[u]=\sup_{\phi\in \Set_n}
\left\{\intom \left(-u\cdot\di \phi-f^*(\phi)\right)dx\right\}.
\]
As $ \left\lbrace \Set_n\right\rbrace_{n=1}^{\infty} $ is a monotonically increasing, $ \Lii n V_{S_n}^f[u]$ exists. Moreover, since $ \Set_n\p \Set_0 $ for all $n\geq 1$, 
\begin{equation}\label{eq: lim 1 V}
\Lii n V_{\Set_n}^f[u]\leq V_{\Set_0}^f[u].
\end{equation}
Let $ \epsilon>0 $ and choose
$ \phi^\epsilon \in \Set_0$ such that 
\[
 V_{\Set_0}^f[u]\leq \epsilon+ \intom \left(-u\cdot\di \phi^\epsilon-f^*(\phi^\epsilon)\right)dx.
\]
Let $ \{\phi^\epsilon_n\} _{n\in \N^*} $ be a sequence converging to $ \phi^\epsilon $ in $ \wo 1 q {\Omega,\R ^{d\times d}} $
and such that $ \phi_n^\epsilon\in \Set_n $ for all $ n\in \N^* $.
Then, 
we use the growth conditions on $  f^* $ in \eqref{eq: growth on f* and grad f*} and \eqref{eq:f^* growth},
the continuity of $f^*$ on its domain and the Lebesgue dominated convergence theorem to obtain that
\[
\intom-f^*(\phi^\epsilon)dx
=\Lii n\intom -f^*(\phi^\epsilon_n)dx.
\]
It follows that
\begin{align*}
V_{\Set_0}^f[u]\leq& \epsilon+ \intom \left(-u\cdot\di \phi^\epsilon-f^*(\phi^\epsilon)\right)dx 
\\
=&
\epsilon+ \Lii n\intom \left(-u\cdot\di \phi^\epsilon_n-f^*(\phi^\epsilon_n)\right)dx
\\
\leq &\epsilon+ \limsup _{n\to \I}V_{\Set_n}^f[u]
\\
= &\epsilon+ \lim _{n\to \I}V_{\Set_n}^f[u].
\end{align*}
As $ \epsilon $ is arbitrary, we have 
\begin{equation}\label{eq: lim 2 V}
\Lii n V_{\Set_n}^f[u]\geq V_{\Set_0}^f[u].
\end{equation}
From  \eqref{eq: lim 1 V} and \eqref{eq: lim 2 V}, we conclude that
 $ \Lii n V_{\Set_n}^f[u]=V_{\Set_0}^f[u]. $
\newline
2.) One has
\begin{align*}
 V_\Set^f[u]=&\sup _{\phi\in\wo 1 q {\Omega,\R^{d\times d} } }
 \left\{\intom \left (-u\cdot \di \phi-f^*(\phi)  \right )dx\right\}
 \\
 =&
\sup _{\phi\in\wo 1 q {\Omega,\R^{d\times d} } }
 \left\{\intom \left (\nabla u\cdot \phi-f^*(\phi)  \right )dx\right\}
 \\
\leq &\intom f(\nabla u)\;dx.
\end{align*}
The inequality above is obtained by using the definition of the Legendre transform $f^*$ of $f$. Let $ \bar \phi\in \pd f (\nabla u) $. Then $ f^*(\bar \phi)+f(\nabla u)=\nabla u\cdot \bar\phi $.
Thanks to the growth conditions \eqref{eq:f growth} and \eqref{eq:pd f growth}
on $ f $, it holds that $ \bar \phi \in
\el q {\Omega,\R ^{d\times d}}$.
Since $\wo 1 q {\Omega,\R ^{d\times d}}$
is dense in $\el q {\Omega,\R ^{d\times d}}$ for the $\el q {\Omega,\R ^{d\times d}}$ norm, we get
\begin{align*}
 \intom f(\nabla u)\;dx=&\intom \left ( 
\nabla u\cdot \bar\phi -f^*(\bar \phi)\right)\;dx
\\
\leq &\sup _{\phi\in\wo 1 q {\Omega,\R^{d\times d} } }
 \left\{\intom \left (\nabla u\cdot \phi-f^*(\phi)  \right )\;dx\right\}
 \\
 =& V^f_\Set[u].
\end{align*}
We conclude that
$ V^f_S[u]= \intom f(\nabla u)\;dx$.
\newline
3.)
The total variation of $ u\in BV(\Omega,\R^{d\times d}) $ 
is
\begin{equation}\label{eq: TV}
 \|Du\|(\Omega) =\sup \left \{ \intom u\cdot \di \phi\;dx: \phi\in C_c^1({\Omega,\R^{d\times d}}); |\phi|\leq 1  \right \}
\end{equation}
 while, using the Legendre transform of $f(\xi)= |\xi|$, we obtain
\begin{equation}\label{eq: Vf}
V^f_S(u)=\sup \left \{ \intom u\cdot \di \phi\;dx: \phi\in \wo1q{\Omega,\R^{d\times d}}; |\phi|\leq 1  \right \}.
 \end{equation}
 It follows directly from \eqref{eq: TV} and \eqref{eq: Vf} that $ \|Du\|(\Omega)\leq  V^f_\Set[u]$. The converse inequality $ \|Du\|(\Omega)\geq  V^f_\Set[u]$ follows from the density of $C_c^1({\Omega,\R^{d\times d}})$ in\\
 $\wo1q{\Omega,\R^{d\times d}}$  and an argument similar to the one made in the proof of (2) in  the proposition.

\cqfd

\section{ Minimization with general displacements.}
\noindent  We consider the following :
\begin{equation}
\label{eq:primal H}
\inf_{(u,\beta)\in \mathscr U^*_\Set}\left\{ I(u,\beta)= V_\Set^f(u)+\int_\Omega\left(  H(\beta)-F\cdot u \right)dx\right\}.
\end{equation}
This problem will be studied via a dual problem that we will formulate  next.
We assume in this section that Assumption \textbf{(A2)} holds.

\subsection{An auxiliary problem}

For $l, k:\mathbb{R}^d  \rightarrow (-\infty,\infty]$, %
define for $u,v\in \R^d$
\begin{align}\label{eq: transform l}
l^\#(v):=\sup_{u\in \bar\Lambda  ,t>0}\left\{u\cdot v-l( u) t-H(t)\right\}
\end{align}
and
\begin{align}\label{eq: transform k}
\quad \quad k_\#(u):=\sup_{v\in \mathbb{R}^d  ,t>0}\left\{(1/t) \left( u\cdot v-k( v) -H(t)\right) \right\}.
\end{align}
Under Assumption \textbf{(A2)}, it is known that $((l^\#)_\#)^\#=l^\#$ and $((k_\#)^\#)_\#=k_\#$ 
(see for instance Lemma A1 of \cite{Gan_uecsc} ).
Call $\mathcal C$ the set of all functions $(k,l)$
with
$k,\;l:\mathbb{R}^d  \rightarrow \mathbb{R} \cup\{\infty\}$  Borel measurable, finite at least at one point, and satisfying  
$l\equiv \infty$ on $\mathbb{R}^d  \setminus \bar \Lambda$ 
and such that
\begin{equation}
\label{eq: inequality l, k}
k(v) +tl( u) +H(t)\geq u\cdot v\quad \forall u,v\in \mathbb{R}^d  , t>0.
\end{equation}
Call $\mathcal C'$ the set of all functions $(k,l)\in \mathcal C$  
such  that
 $ l=k_\# $  and  $ k=l ^{\#} $.
The set $\mathcal C'$ is nonempty. Indeed, $(\chi_{\bar \Lambda}^\#,(\chi_{\bar \Lambda}^\#)_{\#})\in\mathcal C'$ as $((\chi_{\bar \Lambda}^\#)_{\#}))^{\#}=\chi_{\bar \Lambda}^\#$.\\
Let $\mathscr A$ be the set of $(k,l,\varphi)$  such that  $(k,l)\in \mathcal{C}$
and $\phi\in  \Set$.
Consider the following functional defined on $\mathscr A$:
 $$
 J(k,l,\varphi):=\int_\Omega  k(F+\operatorname{div} \phi)\;dx+\int_\Lambda l\;dy+\intom f^*(\phi)\;dx.
 $$ 
The following problem
will play an important role in this section:
\be\label{eq:dual 2}
\inf \left \{ J(k,l,\varphi): (k,l,\varphi)\in \mathscr A\right \}
.
\ee
The value  in Equation \eqref{eq:dual 2}
is the opposite of the value in the following
equation:
\be\label{eq:dual 2 true}
\sup \left \{ -J(k,l,\varphi): (k,l,\varphi)\in \mathscr A\right \}
.
\ee
Let
$\mathscr A'$ 
denote
the subset of $\mathscr A$ consisting of all $ (k,l,\phi)\in \mathscr A$     that satisfy
  $ (k,l)\in \mathcal C' $.
It 
holds that  
\be\label{eq:dual 2 reg}
\inf \left \{ J(k,l,\varphi): (k,l,\varphi)\in \mathscr A\right \}
=
\inf \left \{ J(k,l,\varphi): (k,l,\varphi)\in \mathscr A'\right \}.
\ee 

Indeed, the key observation to this end is that
for $(k,l,\varphi)\in \mathscr A$, one has
$l\geq k_\#$ and $k\geq (k_\#)^\#$ so that
$$
J(k,l,\varphi)\geq J((k_\#)^\#,
k_\#,\varphi)\quad\hbox{and} \quad((k_\#)^\#, k_\#,\varphi)\in \mathscr A'.
$$
For $R>0$, we set 
\[
\mathscr A_R=\left \{ (k,l,\varphi)\in \mathscr A':J(k,l,\varphi)\leq R \right \}.
\]
\begin{lemme}\label{lem:inf_l_bound}
 Assume  $\mathbf{(A1)}$, $\mathbf{(A2)}$ and $\mathbf{(A3)}$ hold. Let $(k,l,\varphi)\in \mathscr A_R$. Set  
$s_l:=-\inf_{u\in\bar\Lambda}l(u)$. Then,
 \begin{equation*}
 \int_{  \Omega}k(F+ \operatorname{div}\varphi )\;dx\geq \mathcal{L}^d{(\Omega)}H^*(s_l)
 -r^* \|F\|_{  L^{1}( \Omega)  }.\label{eq:boundedkn1}
 \end{equation*}
 Moreover, there exists $ M:= M(R,F,f,\Omega,\Lambda)>0 $  such that
 \begin{equation}
\label{eq:boundedsl1}
|s_l|\leq M.
 \end{equation}

\end{lemme}
{Proof.} 
As $\Lambda$ is bounded and $l$ is convex, we choose  $u_l\in\overline{\Lambda}$  such that  $-l(u_l)=s_l$. Since   $k:=l ^{\#}$,  in view of (\ref{eq: transform l}), we have
\begin{equation}\label{eq: ineq H, sl, k}
 -tl(u_l) -H(t)+u_l\cdot v= ts_l -H(t)+u_l\cdot v \leq H^*(s_l)+u_l\cdot v \leq k(v).
 \end{equation}
Using the last inequality in 
{\eqref{eq: ineq H, sl, k}}, one gets
 \begin{align}\label{eq : control int k below00}
  \int_{  \Omega}k(F+ \operatorname{div}\varphi )\;dx
  \geq &\intom \left(H^*(s_l)+u_l\cdot(F+\operatorname{div}\varphi)\right)\;dx
 \\
 = &H^*(s_l)\mathcal{L}^d{(\Omega)}+\intom u_l\cdot F\;dx.\label{eq : control int k below001}
  \end{align}
 We have used the fact that  $u_l$ is a constant vector and $\phi\in \wo 1 q {\Omega,\R^{d\times d} }$ to obtain the equality in 
{\eqref{eq : control int k below001}}.  Hence,
$$  \int_{  \Omega}k(F+ \operatorname{div}\varphi )\;dx\geq \mathcal{L}^d{(\Omega)}H^*(s_l)
  -r^* \|F\|_{L^{1}(\Omega)}.%
   $$
Thus,
\[
R\geq J(k,l,\varphi)\geq -s_l \mathcal{L}^d(\Lambda)+
\mathcal{L}^d{(\Omega)}H^*(s_l)
  -r^* \|F\|_{  L^{1}( \Omega)  }+ \inf f^*.
\]
Thanks to  Lemma \ref{lem:pteOf_h} (iii), $ s_l $  is bounded uniformly in $ l $. 
\begin{flushright}$\square$\end{flushright}

\begin{lemme} 
 Assume  $\mathbf{(A1)}$, $\mathbf{(A2)}$ and $\mathbf{(A3)}$ hold.
\begin{enumerate}
 \item There exists  $ M>0 $  such that for all $ (k,l,\varphi)\in \mathscr A_R$ one has
 \begin{equation}
 \label{eq: L1 bound on l}
  \int_\Lambda |l(y)|dy \leq M.
 \end{equation}

\item There exist $ a_0,b_0,c_0>0 $  such that  for all $ (k,l,\varphi)\in \mathscr A_R$,  the map $ k $ is $ r^* $-Lipschitz, and one has for all $ v\in \mathbb{R}^d $
\be\label{eq:k growth}
-c_0+a_0|v|\leq k(v)\leq b_0+r^*|v|.
\ee
   
\end{enumerate} 

\label{lem:em}
\end{lemme}
Proof. 
1.) 
Recall that for $(k,l,\phi)\in \mathcal A_R$, one has
\[ 
{
 J(k,l,\phi)=\int_\Omega  k(F+\operatorname{div} \phi)\;dx+\int_\Lambda l\;dy+\intom f^*(\phi)\;dx.
 }
\]
By Lemma \ref{lem:inf_l_bound}, for all $(k,l,\phi)\in \mathcal A_R$, if we define $s_l:=-\inf_{u\in\bar\Lambda}l(u)$, we get
\[
R\geq J(k,l,\varphi)\geq  
\mathcal{L}^d{(\Omega)}H^*(s_l)-r^* \|F\|_{  L^{1}( \Omega)  } +\int_{\Lambda} l(y) dy
+ \mathcal{L}^d{(\Omega)}\inf f^*.
\]
Rearranging the terms, we get:
\[
\int_{\Lambda} l(y) dy\leq R-\mathcal{L}^d{(\Omega)}H^*(s_l)+r^* \|F\|_{  L^{1}( \Omega)  } 
- \inf f^*\mathcal{L}^d{(\Omega)}.
\]
By definition of $s_l$ we also have $-s_l\mathcal{L}^d{(\Omega)}  \leq \int_{\Lambda} l(y) dy$
and thus
\begin{equation}\label{eq: control int l}
-s_l\mathcal{L}^d{(\Omega)}  \leq \int_{\Lambda} l(y) dy\leq R-\mathcal{L}^d{(\Omega)}H^*(s_l)+r^* \|F\|_{  L^{1}( \Omega)  } 
- \inf f^*\mathcal{L}^d{(\Omega)}.
\end{equation}
 Lemma \ref{lem:inf_l_bound} also ensures that.
We consider the negative part of $l$ defined by 
$l^-:=\max\{-l,0\}$  and note that
\begin{equation}\label{eq: control int |l|}
 \int_{\Lambda} |l(y)| dy
 =\int_{\Lambda} l(y) dy+2\int_{\Lambda} l^-(y) dy.
 \end{equation}
Observe that, by the definition of $s_l$, we have $l^{-}\leq |s_l|$. This, combined with \eqref{eq: control int l}, \eqref{eq: control int |l|} and  \eqref{eq:boundedsl1} yields \eqref{eq: L1 bound on l}.

\noindent 
2.) Let 
$(k,l,\phi)\in \mathscr A_R$. Since $k=l^\#$,
by equation \eqref{eq: transform l}, $k$ is a $r^*$-Lipschitz
as $\Lambda$ has diameter less or equal to $r^*$.
Next 
\begin{align*}
k(0)=&\sup _{u\in \bar \Lambda, t>0}\left\{ -tl(u)-H(t) \right\}
\\
=&\sup _{ t>0}\left\{ -ts_l-H(t) \right\}.
\end{align*}
As $s_l$ is uniformly bounded, the growth condition
on $H$ ensures that $|k(0)|$ is uniformly bounded say by some 
$b_0>0$. We get then the inequality  
$k(v)\leq b_0+r^*|v|$ for all $v\in \R^d$.
\newline Because of the hypothesis on the domain $\Lambda$,  we take $a_0>0$ such that $B(0,a_0)\p \Lambda$. 
As $(k,l,\phi)\in \mathcal A_R$, we use relation
\eqref{eq: inequality l, k}
to obtain for $v\neq 0$
\begin{align}\label{eq:temp for k}
k(v)\geq v\cdot \left(a_0\frac v{|v|}\right)-l\left(a_0\frac v{|v|}\right)-H(1).
\end{align} 
Thanks to 
equation \eqref{eq: L1 bound on l}, $\int_\Lambda |l|dy$ is
 uniformly bounded  in 
 $ l $. 
 We use in addition the
 fact that $ l $ 
is bounded to 
deduce 
that 
$
     \sup_{y\in \bar B(0,a_0)}|l|(y)
$
is bounded by a constant independent of $l$ (see for instance Theorem 1, p. 236 in
\cite{mtfpf}). Thus equation 
\eqref{eq:temp for k} implies that there exists $c_0>0$
such that 
$k(v)\geq a_0|v|-c_0$ for all $v\in \R^d$.
\cqfd

\begin{pro}\label{lem: existence compressible}
Assume $\mathbf{(A1)}$, $\mathbf{(A2)}$, and $\mathbf{(A3)}$ hold. Assume $\Set$ is a finite dimensional non-empty closed and convex  subset of $ W_0^{1,q }({\Omega,\mathbb{R}^{d\times d}})$.
Then, the functional $ J $ admits a minimizer  $ ( k_0, l_0,\varphi_0)$ in $ \mathscr A' $.
\end{pro}
Proof.
Let $ (\bar k,\bar l,\bar\varphi)\in\mathscr A $. Set $ R=J(\bar k,\bar l,\bar\varphi)$. Take a minimizing sequence\\
$ \{(k_n,l_n,\varphi_n)\}_{n\in \mathbb{N}  ^*} $  of Problem \eqref{eq:dual 2} that is in $ \mathscr A_R $. By Lemma \ref{lem:inf_l_bound} and the growth condition on $f^*$ we may assume without  
loss 
of generality that
$  \{ \varphi_{ n}\}^{\infty } _{ n=1}    $ converges to some $ \varphi_0\in \Set $ weakly in $\el q{\Omega,\mat dd}$. 
Since $\text{Span}(\Set)$ is finite dimensional, $  \{ \varphi_{ n}\}^{\infty } _{ n=1}    $ converges to some $ \varphi_0\in \Set $ strongly in the $\el q{\Omega,\mat dd}$ norm. We deduce
\be\label{eq:intom f^*}
\intom f^*(\phi_0)\;dx\leq \liminf_{n\to \I} \intom f^*(\phi_n)\;dx.
\ee
From Lemma \ref{lem:em}, as $l_n$ is convex, we use Ascoli-Arz\'ela Theorem together with Theorem 1, p. 236 in
\cite{mtfpf} to deduce that up  to a subsequence, we may assume that $ (k_n,l_n) $ converges locally uniformly $\R^d \times \Lambda$ to $ (k_0,l_0)\in \mathcal C' $ on .
The Lebesgue dominated convergence together with 
inequality
\eqref{eq:k growth} yield
\be\label{eq:intom k F phi}
\intom k(F+\di \phi_0)\;dx\leq \liminf_{n\to \I}\intom k_n(F+\di \phi_n)\;dx.
\ee
Since $\seq l n 1$ is uniformly bounded below (thanks to Lemma 
\ref{lem:inf_l_bound} ), by Fatou's Lemma we get 
\be\label{eq:intl}
\int_\Lambda l_0\;dy\leq \liminf_{n\to \I}\int_\Lambda l_n\;dy.
\ee
By inequalities \eqref{eq:intom f^*}, \eqref{eq:intom k F phi}
and \eqref{eq:intl}, we get 
\[
J(k_0,l_0,\phi_0)\leq \liminf _{n\to \infty}J(k_n,l_n,\varphi_n)
\]
and $(k_0,l_0,\phi_0)$ is a minimizer of $J$ over $\mathscr A'$.
\cqfd
\subsection{A uniqueness result}
Here, we prove  the main result of this section. We will need the 
following lemma which is in the spirit of Lemma 4.3 and 
Lemma 4.4 in
\cite{awi-gan}. A proof of Lemma \ref{lem:lcb} is given in 
subsection \ref{sec:app1}.

\begin{lemme}\label{lem:lcb}
Assume assumption \textbf{(A2)} holds.
Consider a lower semicontinuous  function 
$l_0:\mathbb{R}^d  \to \bar{\mathbb{R}}  $
   such that  $\inf_{\bar\Lambda}l_0>-\infty$; $l_0$ is finite on  $\Lambda$ and $l_0\equiv +\infty$ on $\mathbb{R}^d  \setminus\bar\Lambda$. Set $k_0= ({l_0})^{\#}$.
Let $v\in \mathbb{R}^d  $  be such that    $ k_0$ is differentiable at $v$.
\begin{enumerate}
\item \label{lcbi}There exist unique $u_0\in\bar\Lambda$ and $t_0>0$  such that  
$ k_0( v ) =-t_0l_0(u_0)-H(t_0)-u_0\cdot v.$
In addition, $ u_0 $ and $ t_0 $ are characterized by   $u_0=\nabla k_0( v ) $ and $H'(t_0)+l(u_0)=0$.
 \item \label{lcbii} Let $\hat l\in C_b(\mathbb{R}^d  )$ and let $1\geq \epsilon>0$. Define $l_\epsilon=l_0+\epsilon \hat l$ and $k_\epsilon={\left(l_\epsilon  \right)}^{\#}$.
 \begin{enumerate}
  \item \label{lcbiia}There exists a constant $M$ independent of $v$ and $\epsilon$ 
  such that
 \[
  \left|\frac{k_\epsilon(v)-k_0( v ) }{\epsilon}\right|\leq M.
 \]
  \item \label{lcbiib}We have \[ \lim _{\epsilon\rightarrow 0 }   \frac{k_\epsilon(v)-k_0( v ) }{\epsilon}=- t_0\hat l(u_0).\]
 \end{enumerate}

\end{enumerate}

\end{lemme} 

\noindent 
Next, we give the main result of this section.
\begin{thm}\label{thm:main 1}
Assume $\mathbf{(A1)}$, $\mathbf{(A2)}$, and $\mathbf{(A3)}$ hold. Assume $\Set$ is a finite dimensional non-empty closed and convex  subset of $ W_0^{1,q }({\Omega,\mathbb{R}^{d\times d}})$. Assume $F$ satisfies the condition  $\textbf{(ND)}_{\Set}$. Then,
Problems \eqref{eq:primal H}  
and \eqref{eq:dual 2 true} are dual.
Problem  \eqref{eq:dual 2 true} admits a maximizer $ (k_0,l_0,\varphi_0) $ with $ k_0=l_0^\# $   and $ l_0=(k_{0})_\# $.
Problem  \eqref{eq:primal H} admits a
unique minimizer  $ (u_0,\beta_0) $. Moreover $ u_0 $ satisfies
\begin{align*}
\begin{cases}
u_0=&\nabla k_0(F+\operatorname{div} \varphi_0)\\
\varphi_0\in&\Phi_S(u_0).
\end{cases}
\end{align*}

\end{thm}
{Proof.} 
{\bf Step 1}.
For 
$(u,\beta)\in \mathscr U^*_\Set$  and $(k,l,\varphi)\in \mathscr A$, one has
\begin{align*}
  I(u, \beta)
 =
  &V_\Set^f(u)+\intom\left( H(\beta)-F\cdot u \right)\;dx
  \\
  \geq&  
  \intom \left( -u\cdot (\mathop{div} \phi+F)\right)dx-\intom f^*(\phi)\;dx
  \\&+
  \intom H(\beta)dx+\intom\beta l (u)\;dx 
  -\int_\Lambda l(y)\;dy
  \\
  \geq&  
  \intom -k (\mathop{div} \phi+F)\;dx
  -\intom f^*(\phi)\;dx -\int_\Lambda l(y)\;dy.
\end{align*}
Thus $I(u,\beta)\geq -J(k,l,\phi)$ with equality
 if and only if  
 {$\varphi\in \Phi_\Set(u)$} and
\[
 k(F+\operatorname{div} \varphi)+\beta l(u)+H(\beta)=u\cdot(F+\operatorname{div} \varphi).
\]
Note that if  $k$ is convex, the map $ \nabla k (F+\operatorname{div} \varphi)$ is well defined as the map $ F+\operatorname{div} \varphi $ is non-degenerate.
Using Lemma \ref{lem:lcb}(i), it follows that if $k$ is convex, then 
$I(u,\beta)=-J(k,l,\phi)$  if and only if
\be\label{eq:I=J}
\begin{cases}
\phi &\in \Phi_\Set(u)
\\
 u&=\nabla k(F+\operatorname{div} \varphi)
 \\
 \beta&=(H')\s(-l(u))
\end{cases}.
\ee
{\bf Step 2}.
Thanks to equation
\eqref{eq:dual 2 reg},
we may find
a  maximizer $ (k_0,l_0,\varphi_0) $ of Problem  \eqref{eq:dual 2}
satisfying $ k_0=l_0^\# $   and $ l_0=(k_0)_{\#} $.
The 
function 
$u_0=\nabla k_0(F+\operatorname{div} \varphi_0)$ is well defined 
as $k_0$ is convex and we set   $\beta_0 =(H')\s(-l(u_0))$. We are to show that 
$(u_0,\beta_0)\  \in \mathscr U^*_\Set$ and $\phi_0  \in \Phi_\Set(u_0)$.
\\
{\bf Step 3}.
Let $ \bar l\in C_c(\mathbb{R}^d  ) $. For $ \epsilon\in (0,1) $,  define $ l_\epsilon=l_0+\epsilon \bar l $  and $ k_\epsilon=(l_\epsilon)^\# $. Using  Lemma \ref{lem:lcb}, one has
\begin{equation} \label{eq:lim_intom_k_epsilon_bis}
\begin{split}
 &\lim _{\epsilon\rightarrow 0^+} \int_\Omega  (1/\epsilon)\left( k_0(F+\operatorname{div} \phi_0)-k_\epsilon(F+\operatorname{div} \phi_0)\right) dx\\
=&\int_\Omega  \beta_0\bar l(\nabla k_0(F+\operatorname{div} \phi_0))\;dx=\int_\Omega  \beta_0\bar l(u_0)\;dx.
\end{split}
\end{equation}
Since $ J(k_0,l_0,\varphi_0)\leq J(k_\epsilon,l_\epsilon,\varphi_0) $,
we deduce that 
$ -\int_\Lambda \bar l\;dy +\int_\Omega  \beta_0\bar l(u_0)\;dx\leq 0$. As we can replace $ \bar l $ by $ -\bar l $, one deduces
that
$ \int_\Lambda \bar l \;dy=\int_\Omega  \beta_0\bar l(u_0)\;dx$. 
Therefore $ (u_0,\beta_0)\in \mathscr U_S^* $.
\\
{\bf Step 4}.
Let   $ \varphi\in \Set $. For $ \epsilon\in (0,1) $,  set   $ \varphi_\epsilon=\epsilon \varphi+(1-\epsilon) \varphi_0 $.  
By the convexity of $ \Set $, the map $ \varphi_\epsilon $ belongs to $\Set$.
As $ J(k_0,l_0,\varphi_0)\leq J(k_0,l_0,\varphi_\epsilon) $, we have
\be\label{eq:ineq1}
\begin{aligned}
0\leq&\int_\Omega    (1/\epsilon)\left(
k_0(F+\operatorname{div} \varphi_0+\epsilon\operatorname{div}(\varphi-\varphi_0))
-
k_0(F+\operatorname{div} \varphi_0)\right)\;dx
\\
&+(1/\epsilon) \int_\Omega\left(f^*(\phi_0+\epsilon(\phi-\phi_0))
-f^*(\phi_0)\right)\;dx
\end{aligned}
\ee 
Thanks to Lemma \ref{lem:lcb}, Inequality \eqref{eq:ineq1} implies 
\begin{align*}
 &\int_\Omega \left( u_0\cdot \operatorname{div}(\varphi- \varphi_0)
+grad[ f^*](\phi_0)\cdot (\varphi- \varphi_0)\right)\;dx\\
=&\int_\Omega\left(  \nabla k_0(F+\operatorname{div} \varphi_0)\cdot\operatorname{div}(\varphi- \varphi_0) +grad[ f^*](\phi_0)\cdot (\varphi- \varphi_0)\right)\;dx
\\
\geq& 0.
\end{align*}
It follows from Proposition \ref{pro:convex pseudo projected} that $\varphi_0\in \Phi_\Set(u_0) $.
\\
{\bf Step 5}. Since $(u_0,\beta_0)\  \in \mathscr U^*_\Set$, $\phi_0  \in \Phi_\Set(u_0)$,  $u_0=\nabla k_0(F+\operatorname{div} \varphi_0)$, and $\beta_0 =(H')\s(-l(u_0))$,
we deduce that $I(u_0,\beta_0)=J(k_0,l_0,\phi_0)$ and 
$u_0$ is a minimizer of Problem \eqref{eq:primal H} thanks to relation \eqref{eq:I=J}. Suppose $(u_1,\beta_1)\  \in \mathscr U^*_S$ is another minimizer
of Problem \eqref{eq:primal H}. Then we have $I(u_1,\beta_1)=J(k_0,l_0,\phi_0)$ 
and by relation \eqref{eq:I=J}, we get $u_1=\nabla k_0(F+\operatorname{div} \varphi_0)$ which implies
$u_1=u_0$. Next the strict convexity of $H$ yields that $\beta_0=\beta_1$. We conclude that $(u_0,\beta_0)$ is the unique minimizer of Problem \eqref{eq:primal H}
and $u_0$ is characterized by
\begin{align*}
\begin{cases}
u_0=&\nabla k_0(F+\operatorname{div} \varphi_0)\\
\varphi_0\in&\Phi_\Set(u_0).
\end{cases}
\end{align*}
\begin{flushright}$\square$\end{flushright}

\begin{cor}\label{cor:1}
Assume $\mathbf{(A0)}$, $\mathbf{(A1)}$, $\mathbf{(A2)}$, and $\mathbf{(A3)}$ hold. Assume $\Set$ is a finite dimensional non-empty closed and convex  subset of $ W_0^{1,q }({\Omega,\mathbb{R}^{d\times d}})$  and $\nabla\phi$ is non-degenerate whenever $\phi\in \Set$. Suppose $ F $ has a countable range (thus degenerate) Then,  $F$ satisfies the condition  $\textbf{(ND)}_{\Set}$  and   
problem \eqref{eq:primal H} admits a unique solution.
\end{cor}

\begin{cor}\label{cor:2}
Assume $\mathbf{(A1)}$, $\mathbf{(A2)}$, and $\mathbf{(A3)}$ hold. Assume $\Set$ is a finite dimensional subspace of $ W_0^{1,q }({\Omega,\mathbb{R}^{d\times d}})$ and $\nabla \phi$ has a countable range whenever $\phi\in \Set$. Suppose $ F $ is non-degenerate. Then,  $F$ satisfies the condition  $\textbf{(ND)}_{\Set}$  and   
problem \eqref{eq:primal H} admits a unique solution.
\end{cor}

\section{ The incompressible case }
\noindent Throughout this section, we assume that $\Set$ is a  subset of  $\wo1q{\Omega,\R^{d\times d}}$.  We consider the following  problem:
\be\label{eq:primal incompressible}
\inf_{u\in \mathcal{U}_\Set^1}\left\{ I_0(u):=V^f_\Set( u)-\intom F\cdot u\; dx \right\}
\ee
and we recall that the set $\mathcal{U}_\Set^1$ is defined as
\[
 \mathscr U_\Set^1=\left\{u\in \mathscr U_\Set: \int_\Omega l(u(x))\;dx=\int_\Lambda l(y)\;dy \; \forall l \in C_c(\mathbb{R}  ^d)\right\}.
\]
We assume  $\lebd(\Omega)=\lebd(\Lambda)$ so that $ \mathscr U_\Set^1$ is non-empty.

\subsection{Existence and uniqueness  via duality}
We  study Problem
\eqref{eq:primal incompressible} via duality.
Let $u\in \mathscr U^1_\Set$, $\varphi \in \Set$,  $l\in C(  \Lambda )$ and $k:\mathbb{R}^d  \to \mathbb{R}  $
satisfy $k(v)+l(u)\geq u\cdot v$ for all $u\in \Lambda$ and all $v\in \mathbb{R}^d  $.
One has
\begin{align}
  &V^f_\Set (u)-\int_\Omega  F\cdot u\; dx
  \\
  =& -\int_\Omega  u\cdot(F+\operatorname{div} \varphi)\;dx+\intom l(u)\;dx-\int_\Lambda l(y)\;dy-\intom f^*(\varphi)\;dx\label{eq:ineq dual comp 1}
 \\
 \geq & -\int_\Omega  k(F+\operatorname{div} \varphi)\;dx-\int_\Lambda l(y)dy
 -\intom f^*(\varphi)\;dx.\label{eq:ineq dual comp 2}
\end{align}
This suggests that we consider the dual problem
\begin{equation}\label{eq:pro dual gro 1 H 1}
M_0:=\inf_{( k,l,\varphi)\in A_0} \left \{ J(k,l,\varphi):=\int_\Omega  k(F+\operatorname{div} \varphi)\;dx+\int_\Lambda l(y)\;dy +\intom f^*(\varphi)\;dx \right \}
\end{equation}
with $ A_0 $ being the set of all $ (k,l,\phi) $ such that $ \varphi\in \Set $, $l\in C(  \Lambda  )$, $ \inf_\Lambda l=0 $ and  $k:\mathbb{R}^d  \to \mathbb{R}  $
satisfies $k(v)+l(u)\geq u\cdot v$ for all $u\in \bar\Lambda , $ and all $v\in \mathbb{R}^d  $. 
Remark that we have 
\begin{equation}\label{eq:pro dual gro 1 H 1 true}
-M_0=\sup_{( k,l,\varphi)\in A_0} \left \{ -J(k,l,\varphi) \right \}.
\end{equation}

\subsubsection{Existence and regularity of minimizers of Problem {\eqref{eq:pro dual gro 1 H 1}}}

\newcommand{\ds}{\mathcal{C}}

Denote by  $ \ds $ the set 
of all $ (k,l) $  such that 
$ k:\R^d\to \R $ and $ l:\R^d\to \R\cup\{\infty\} $ satisfy
\be\label{eq:ineq}
k(v)+l(u)\geq u\cdot v;\quad \forall u\in \bar\Lambda;\quad \forall v\in \R^d
\ee
and $l\equiv \infty$ on $\R^d\setminus \bar \Lambda$.
Consider the subset $\ds_0$ of $\ds$ consisting of  $ (k,l)\in \ds $  such that $l\in C(  \Lambda  )$ and $ \inf_\Lambda l=0$.
The following lemma is standard:

\begin{lemme}\label{lem:reg incompressible}
Let  $(k,l)\in \ds$. It holds that $(l^*, l^{**})\in \ds$,
 $l^*\leq k$, $0\leq l^{**}\leq l$ and $l^{***}=l^*$. 
If $(k,l)\in \ds_0$ then $l^*(0)=0$.
\end{lemme}
Let us denote by $C_0' $ the set of all $(k,l)\in C_0$
such that $l^*=k$, $k^*=l$, $k(0)=0$, and $l\geq 0$, and by $A'_0$ 
the set of all $(k,l,\phi)$ with $(k,l)\in C_0'$
and $\phi\in \Set$.
Remark that an element in $C_0' $ is the couple $(\chi_{\bar\Lambda},(\chi_{\bar\Lambda})^*)$. Hence 
$A'_0$ is nonempty when $\Set$ is nonempty.
One readily checks that, in light of Lemma \ref{lem:reg incompressible},
Problem \eqref{eq:pro dual gro 1 H 1}
has the same infimum value as
\begin{equation}\label{eq:pro dual gro 1 H 1 reg}
 \inf_{(k,l,\phi)\in A_0'} \left \{J(k,l,\varphi):= \int_\Omega  k(F+\operatorname{div} \varphi)\;dx+\intom f^*(\varphi)\;dx+\int_\Lambda l(y)\;dy \right \}.
\end{equation}
We  recall that
$r^*$ is  such 
 that $B(0,1/r^*)\p \Lambda\subset   B(0,r^*/2)$; 
\begin{lemme}\label{lem:existence J incompressible}
Assume  $\mathbf{(A1)}$ and $\mathbf{(A3)}$ hold. Assume that the set $\Set$ is a finite dimensional non-empty closed and convex  subset of  $ W_0^{1,q }({\Omega,\mathbb{R}^{d\times d}})$. Then, problem
\eqref{eq:pro dual gro 1 H 1 reg} admits a minimizer   $ (k_0,l_0,\varphi_0)\in A_0' $ 
with $ k_0 $ convex and $ r^* $-Lipschitz and $ k_0(0)=0 $. %
\end{lemme}
{Proof.}
Consider a minimizing sequence 
$\left\lbrace (k_n,l_n,\phi_n)\right\rbrace_{n=1}^{\infty}$
of Problem \eqref{eq:pro dual gro 1 H 1 reg}.
Since $k_n=l_n^*$ and $l_n=(k_n)^*$, $k_n$ is $r^*$-Lipschitz. As $k_n(0)=0$, we use Ascoli-Arz\'ela theorem to deduce that a subsequence of
$\left\lbrace k_n\right\rbrace_{n=1}^{\infty}$ converges locally uniformly to some $ k_0$.
Next, using the growth condition \eqref{eq:f^* growth} on $f^*$ as well as the facts that $k_n$
is $r^*$-Lipschitz, $k_n(0)=0$, we establish the following estimate :
\begin{equation}\label{eq : estimate from dual Pro}
J(k_n,l_n,\phi_n)\geq \intom \left(-r^*|F+\di \phi_n|+c^p\frac{|\phi_n|^q}{q} -b\right)\;dx+\int_\Lambda l_n(y)\;dy
\end{equation}

As the left hand side of (\ref{eq : estimate from dual Pro}) is bounded, $l_n\geq 0$ and  $\Set$ is finite dimensional,
we deduce from \eqref{eq : estimate from dual Pro} that a subsequence of $\left\lbrace \phi_n\right\rbrace_{n=1}^{\infty}$ converges strongly to some
$\phi_0$ in  $\wo1q{\Omega,\R^{d\times d}}$. Invoking \eqref{eq : estimate from dual Pro} again, we show that $\left\lbrace \int_\Lambda l_n(y)dy\right\rbrace_{n=1}^{\infty} $ is bounded. This, combined with  the fact that $l_n$ is  non negative and convex, yields the existence of a subsequence of
$\left\lbrace l_n\right\rbrace_{n=1}^{\infty}$ that converges locally uniformly to some $ l_0$ (see for instance Theorem 1, p. 236 in
\cite{mtfpf}).
One readily checks that $(k_0, l_0, \phi_0)\in A_0'$.
We next exploit lower semi-continuity 
properties of the functional $J$ to conclude that
$(k_0, l_0, \phi_0)$ is a minimizer of $J$ over $A_0'$.
\begin{flushright}$\square$\end{flushright}

\subsubsection{A duality result}
We have the following theorem.
\begin{thm}\label{thm:main 2}
Assume $\mathbf{(A1)}$ and $\mathbf{(A3)}$ hold. Assume $\Set$ is a finite dimensional non-empty closed and convex  subset of   $  W_0^{1,q }({\Omega,\mathbb{R}^{d\times d}})$. Suppose that the map $ F $ satisfies the condition $\textbf{(ND)}_{\Set}$.
Then
Problems \eqref{eq:primal incompressible}  
and \eqref{eq:pro dual gro 1 H 1 true} are dual.
Problem  \eqref{eq:pro dual gro 1 H 1 true} admits a maximizer $ (k_0,l_0,\varphi_0) $ with $ k_0=l_0^* $   and $ l_0=(k_0)^* $.
Problem  \eqref{eq:primal incompressible} admits a
unique minimizer  $ u_0 $. Moreover $ u_0 $ satisfies
\begin{align*}
\begin{cases}
u_0=&\nabla k_0(F+\operatorname{div} \varphi)\\
\varphi_0\in&\Phi_S(u_0)
\end{cases}.
\end{align*}
\end{thm}
Proof.
Suppose $ u\in \mathscr U_S^1 $ and  $( k,l,\varphi)\in A_0 $.
Using 
\eqref{eq:ineq dual comp 1}
and 
\eqref{eq:ineq dual comp 2},
we see that 
$I_0(u)\geq -J(k,l,\varphi)$
with equality
 if and only if 
$ \varphi\in \Phi_S(u) $ and 
 $l(u)+k(F+\di \phi)=u\cdot (F+\di \phi) $ for almost every $ x\in \Omega $. 
The latter condition reduces to
 $ u(x)=\nabla k(F(x)+\operatorname{div} \varphi(x)) $ if $k$ is convex, under the assumption  $F+\operatorname{div} \varphi$ is non-degenerate. Now, let $ (k_0,l_0,\varphi_0) \in A_0'$ be a minimizer of $ J $ over $ A_0 $.
Since $ F+\operatorname{div} \varphi_0 $ is non-degenerate and $ k_0 $ is convex, the map $ u_0=\nabla k_0(F+\operatorname{div}\varphi_0) $ is well defined.
\newline
\newline
{\bf Variation around {$l_0$}.}
Let $ \bar l\in C_c(\mathbb{R}^d  ) $. For $ \epsilon\in (0,1) $,  set $ l_\epsilon=l_0+\epsilon \bar l $ and $ k_\epsilon=(l_\epsilon)^*. $ Let $ v\in \mathbb{R}  ^d $ be a point where $ k_0 $ is differentiable. 
Using the measurable selection theorem, one deduces that there
exists 
$ T_\epsilon:\mathbb{R}  ^d\to \mathbb{R}^d   $ measurable  such that 
for all  $ \epsilon\in [0,1) $ 
\[
k_\epsilon
(v)=T_\epsilon(v)\cdot v-l_\epsilon(T_\epsilon(v)).\]
Then, for $ \epsilon\in (0,1),$ we have
\begin{align}\label{eq: var 1}
\bar l(T_\epsilon(v))\leq-(1/\epsilon)\left( k_\epsilon(v)-k_0( v )\right) 
\leq\bar l(T_0(v))
\end{align}

and
\begin{align}\label{eq: var 2}
\left|(1/\epsilon)\left( k_\epsilon(v)-k_0( v )\right)\right|
\leq\|\bar l\|_{  L^{\infty}(\mathbb{R}^d  )  }.
\end{align}
Moreover,
\begin{align}\label{eq: var 3}
 \lim _{\epsilon\rightarrow {0^+} } -(1/\epsilon)\left( k_\epsilon(v)-k_0( v )\right)
=\bar l(T_0(v)).
\end{align}
We refer the reader to Lemma \ref{lem:lcb:H0} 
for \eqref{eq: var 1}- \eqref{eq: var 3}. Hence, as
$$ T_0(F+\operatorname{div} \psi_0)=\nabla k_0(F+\operatorname{div} \psi_0)=u_0\quad a.e.$$ 
using again \eqref{eq: var 3},
one has
 \begin{equation}
 \begin{aligned} \label{eq:lim_intom_k_epsilon_bis_bis}
 &\lim _{\epsilon\rightarrow 0^+} \int_\Omega (1/{\epsilon})\left( k_0(F+\operatorname{div} \psi_0)-k_\epsilon(F+\operatorname{div} \psi_0)\right) dx\\
=&\int_\Omega  \bar l(T_0(F+\operatorname{div} \psi_0))\;dx=\int_\Omega  \bar l(u_0)\;dx.
\end{aligned}
\end{equation}
Since $ J(k_0,l_0,\varphi_0)\leq J(k_\epsilon,l_\epsilon,\varphi_0) $,
we deduce from \eqref{eq:lim_intom_k_epsilon_bis_bis} that\newline 
$ -\int_\Lambda \bar l +\int_\Omega  \bar l(u_0)\leq 0$.
By replacing $ l $ by $ -l $ in the above argument,
one deduces
that
$ \int_\Lambda \bar l =\int_\Omega  \bar l(u_0)$. 
As a result, $ u_0\in \mathscr U_S^1 $.
\newline\newline
{\bf Variation around $\varphi_0$.}
Let $ \varphi\in \Set $. 
For $ \epsilon\in (0,1) $, by convexity of $ \Set $, we have
$ \varphi_\epsilon:=\epsilon \varphi+(1-\epsilon) \varphi_0\in \Set $.  
Then $ J(k_0,l_0,\varphi_0)\leq J(k_0,l_0,\varphi_\epsilon) $. 
This implies that
\begin{align*}
0\geq\int_\Omega (1/\epsilon) \left( k_0(F+\operatorname{div} \varphi_0)-k_0(F+\operatorname{div} \varphi_0+\epsilon\operatorname{div}(\varphi-\varphi_0))\right)&\\ +
f^*(\phi_0)-f^*(\phi_0+\epsilon(\phi-\phi_0)\;dx
.
\end{align*}
As $\epsilon$ tends to $0^+$, the above equation yields 
\begin{align*}
0\geq&-\int_\Omega  \nabla k_0(F+\operatorname{div} \varphi_0)\cdot\operatorname{div}(\varphi- \varphi_0) -grad[ f^*](\phi_0)\cdot (\varphi- \varphi_0)\;dx
\end{align*}
It follows from Proposition \ref{pro:convex pseudo projected} that $\varphi_0\in \Phi_\Set(u_0) $.
\begin{flushright}$\square$\end{flushright}

\begin{cor}\label{cor:21}
 Assume $\mathbf{(A0)}$, $\mathbf{(A1)}$, and $\mathbf{(A3)}$ hold. Assume $\Set$ is a finite dimensional non-empty closed and convex  subset of $ W_0^{1,q }({\Omega,\mathbb{R}^{d\times d}})$  and $\nabla\phi$ is non-degenerate whenever $\phi\in \Set$.
Suppose $ F $ has a countable range (thus degenerate). Then,  $ F $  satisfies the condition $\textbf{(ND)}_{\Set}$ and   
problem \eqref{eq:primal incompressible} admits a unique solution.
\end{cor}

\begin{cor}\label{cor:22}
Assume   
$\mathbf{(A1)}$
and $\mathbf{(A3)}$ hold. Assume $\Set$ is a finite dimensional 
subspace of $ W_0^{1,q }({\Omega,\mathbb{R}^{d\times d}})$  and
$\nabla \phi$ has a countable range whenever $\phi\in \Set$. Suppose $ F $ is non-degenerate. Then, $ F $  satisfies the condition $\textbf{(ND)}_{\Set}$ and   
problem \eqref{eq:primal incompressible} admits a unique solution.
\end{cor}

\subsection{A Link between Problem {\eqref{eq:primal H}} and Problem {\eqref{eq:primal incompressible}}}
Here, we explore  the relationships between 
problem {\eqref{eq:primal H}} and problem {\eqref{eq:primal incompressible}}. For this purpose, we make a further assumption of the domains  $\Omega$ and $\Lambda$ by requiring that $\Omega=\Lambda$. Assume $\mathbf{(A1)}$ holds and recall $\left\lbrace H_n\right\rbrace_{n=0}^{\infty}$  as defined in \eqref{eq H_0} and \eqref{eq H_n}. Then, Lemma \ref{lem: lim H0} ensures that  
$\mathbf{(A2)}$%
holds for $H_n$  for all $n\geq 1$.
Define
\[
 I_n(u,\beta):=V_\Set^f(  u)+\intom H_n(\beta)-u\cdot F\;dx\qquad 
 {n\geq 1}
 \label{eq:primal n}
\]
and  
\[
 I_0(u ):=V_\Set^f(  u)-\intom u\cdot F\;dx.
 \label{eq:primal 0}
\]
Recall  that $C_0$ is the set of all $(k,l)$
such that  $l\in C(\bar \Lambda)$, 
$\inf l=0$
and $k:\R^d\to \R$ satisfies
for all $u\in \Lambda$ and all $v\in \R^d$:
\be\label{eq:C_0}
k(v)+l(u)\geq u\cdot v.
\ee
Let $C_n$ be the set of all $(k,l)$
such that  $l\in C(\bar \Lambda)$ and $k:\R^d\to \R$ 
satisfy:
\be\label{eq:C_n}
k(v)+tl(u)+H_n(t)\geq u\cdot v;\quad \forall u\in \Lambda;\;\forall v \in \R^d.
\ee
We denote  
{by} $\mathcal A_0$ the set of all $(k,l,\phi)$ satisfying $(k,l)\in C_0$
and $\phi\in S$.   Similarly $\mathcal A_n$ denotes 
{the} set of all $(k,l,\phi)$ satisfying $(k,l)\in C_n$
and $\phi\in S$.
If  
{$(k,l,\phi)\in \mathcal{A}_0\cup \mathcal{A}_n$}, we still set 
\[
 J(k,l,\varphi)=\int_\Omega  k(F+\operatorname{div} \varphi)\;dx+\int_\Lambda l(y)\;dy +\intom f^*(\varphi)\;dx.
\]

\bl
Assume $\mathbf{(A1)}$, $\mathbf{(A2)}$, and $\mathbf{(A3)}$ hold. Assume $\Set$ is a finite dimensional non-empty closed and convex  subset of   $  W_0^{1,q }({\Omega,\mathbb{R}^{d\times d}})$.  
For each $n\in \N$, let $(u_n,\beta_n)$ be the unique minimizer of $I_n$ over $\mathscr U_\Set^*$ as given by Theorem \ref{thm:main 1} and  let $(k_n,l_n,\phi_n)$ be a minimizer of $J$ over  
{$\mathcal{A}_n$} with  
{$k_n$  convex}
and $r^*$-Lipschitz as ensured  by Proposition \ref{lem: existence compressible} and Lemma \ref{lem:existence J incompressible}.
Then,
\begin{enumerate}
 \item The sequence $\{I_n(u_n,\beta_n)\}_{n\in \N^*}$
 is bounded.
 \item The sequence $\{ \beta_n \}_{n\in \N^*}$ 
 converges to $1$ in $\el 2 \Omega$.
 \item %
 The sequence $\{ \phi_n \}_{n\in \N^*}$ 
 admits a subsequence  that 
 converges to some $\bar \phi$ in $S$ with respect to the $W_0^{1,q}\left(\Omega,\mathbb{R}^{d\times d} \right)-$ norm.
 \end{enumerate}
\label{lem:various compactness}
\El
Proof.  
{\bf Step 1.}
Let $\bar u\in \mathscr U_\Set^1$. We have 
{$(\bar{u},1)\in \mathscr U_\Set^*$} and thus $I_n(u_n,\beta_n)\leq I_n(\bar u,1)$ for all $n\geq 1$.
As $H_n(1)=0$, it holds that $I_n(\bar u,1)=V_\Set^f(\bar u)-\intom \bar u\cdot F\;dx$ which is finite.
Hence
\begin{equation}\label{eq0 : bound on In}
R_0:=V_\Set^f(\bar u)-\intom \bar u\cdot F\;dx\geq I_n(u_n,\beta_n).
\end{equation}
On the other hand, we use growth condition \eqref{eq:f^* growth} to get
\begin{equation}\label{eq1 : bound on In}
 I_n(u_n,\beta_n)\geq \intom (-b+u_n\cdot F)\; dx\geq -b\lebd(\Omega)-r^*\|F\|_{\el 1 {\Omega,\R^d}} := - R_1.
\end{equation}
Finally, we use \eqref{eq0 : bound on In} and \eqref{eq1 : bound on In} to  prove (1).
\newline
{\bf Step 2.} Let  $\phi_0 \in \Set$. As $u_n$ has values in $\Lambda$, it holds that
\begin{equation}\label{eq2 : bound on In}
 V_\Set^f(u_n)=\sup_{\phi\in S}\intom(-u_n\di \phi-f^*(\phi))\;dx\geq  \intom(-r^*|\di \phi_0|-f^*(\phi_0))\;dx=: R_2.
\end{equation}
and
\begin{equation}\label{eq3 : bound on In}
 \intom -u_n\cdot F\;dx\geq -r^*\|F\|_{\el 1 {\Omega,\R^d}}.
\end{equation}
 We combine \eqref{eq0 : bound on In}, \eqref{eq1 : bound on In}, \eqref{eq2 : bound on In}, \eqref{eq3 : bound on In} to get 
\begin{equation}\label{eq : bound on H_n(bn)}
R_2 -r^*\|F\|_{\el 1 {\Omega,\R^d}}+\intom H_n(\beta_n)\;dx\leq I_n(u_n,\beta_n)\leq R_0.
\end{equation}

 Setting $c_0 \lebd(\Omega):= R_0- R_2+r^*\|F\|_{\el 1 {\Omega,\R^d}},$ we use lemma \ref{lem: lim H0} and  
\eqref{eq : bound on H_n(bn)} to obtain
\begin{align*}
\intom n(\beta_n(x)-1)^2 dx\leq (c_0+\bar H-H(1))\lebd(\Omega).
\end{align*}
This establishes (2).
\newline
{\bf Step 3.} As $\left\lbrace H_n \right\rbrace_{n=1}^{\infty} $ is a non decreasing sequence that converges to $H_0$, it holds that $C_{n+1}\p C_n\p C_0$ for all $n\in \N$. 
Thus, as $(k_n,l_n)\in C_{n}$,  we have $(k_n,l_n)\in C_{0}$ so that
\begin{align}\label{eq : ineq C_0}
k_n(F+\di \phi_n) + l_n(x)\geq   x\cdot(F+\di \phi_n).
\end{align}
 Since $-J(k_n,l_n,\phi_n)=I_n(u_n,\beta_n)$, we have $J(k_n,l_n,\phi_n)\leq R_1$ for all $n\in \N^*$. This, combined with $\Omega=\Lambda$, and \eqref{eq : ineq C_0}  yields 
\begin{align}\label{eq1 : bound phi_n}
 R_1\geq& \intom \left(k_n(F+\di \phi_n)+l_n(x)+f^*(\phi_n)  \right)dx \\
 \geq& \intom \left(x\cdot(F+\di \phi_n)+f^*(\phi_n)  \right)dx.
 \end{align}
 In view of the growth condition \eqref{eq:f^* growth} and boundedness of $\Omega$, \eqref{eq1 : bound phi_n} implies
 \begin{align}
 R_1\geq \intom \left(r^*|F+\di \phi_n|-b+c^p\frac{|\phi_n|^q}{q} \right)dx.
 \end{align}

As, $\Set$ is of finite dimension and the $\di$ operator is continuous on $\Set$, we conclude that $\left\lbrace\phi_n \right\rbrace_{n=1}^{\I} $ is convergent  up to a subsequence in $W_0^{1,q}\left(\Omega,\mathbb{R}^{d\times d} \right) $ which 
allows us to conclude (3).
\cqfd

\bt\label{thm:main 3}
Assume $\mathbf{(A1)}$, $\mathbf{(A2)}$, and $\mathbf{(A3)}$ hold. Assume $\Set$ is a finite dimensional non-empty closed and convex  subset of $ W_0^{1,q }({\Omega,\mathbb{R}^{d\times d}})$.  
Assume  $F$ satisfies the condition $\mathbf{(ND)_{\mathcal{S}}}$. For each $n\in \N$, let $(u_n,\beta_n)$ be the unique minimizer of $I_n$ over $\mathscr U_\Set^*$ as given by Theorem \ref{thm:main 1} and  let $(k_n,l_n,\phi_n)$ be a minimizer of $J$ over 
$\mathcal{A}_n$%
with $k_n$ is convex and $r^*$-Lipschitz as ensured  by Proposition \ref{lem: existence compressible} and Lemma \ref{lem:existence J incompressible}.   
Suppose that $k_n$ is differentiable for all $n\in \N^*$.
Then, the sequence $\{ u_n \}_{n\in \N^*}$ 
 converges almost everywhere to the unique minimizer $u_0$ of
 $I_0$ over 
   $\mathscr U_S^1$. In addition, the minima $\left\lbrace I_n(u_n,\beta_n)\right\rbrace_{n=1}^{\infty} $ 
   converge to $I_0(u_0)$.
   \et
Proof. 
{\bf Step 1.}
For $n\in \N^*$, set  
$\bar k_n=k_n-k_n(0)$. Note that  
we have $\bar k_n(0)=0$. 
Since the  functions $k_n$ are $r^*$-Lipschitz, so are the  
functions $\bar k_n$ and
we obtain that, up to a subsequence, 
the sequence 
$\left\lbrace \bar k_n\right\rbrace_{n=1}^{\infty}$ converges locally uniformly to
a certain  
function 
$\bar k$.
Since $F+\di \phi_n$ is non-degenerate, we have that $\nabla \bar k_n(F+\di \phi_n)$
is well-defined. Furthermore, Lemma \ref{lem:various compactness} ensures that  $\left\lbrace \phi_n\right\rbrace_{n=1}^{\infty} $ converges 
up to a subsequence 
to some $\bar \phi\in S$ with respect to 
the
$W^{1,q}\left(\Omega, \mathbb{R}^d \right)-$ 
{norm}. As a result,   $\left\lbrace \di \phi_n\right\rbrace_{n=1}^{\infty}$ converges to $\di \bar\phi$ in $L^q(\Omega,\mathbb{R}^d)$. Since $\Set$ is of finite dimension, the $L^q$ convergence of  $\left\lbrace \di \phi_n\right\rbrace_{n=1}^{\infty}$  reduces to a pointwise convergence. Next, using the convexity of the 
$\bar k_n$ and the  pointwise convergence of $\left\lbrace \di \phi_n\right\rbrace_{n=1}^{\infty}$  to  
{$\di \phi $}, we deduce that up to  a subsequence $\left\lbrace \nabla \bar k_n(F+\di \phi_n)\right\rbrace_{n=1}^{\infty} $ converges a.e to $ \nabla\bar k(F+\di \bar\phi)$ ( cf. \cite{Roc_ca} Theorem 25.7).\\ 
As a duality result, Theorem \ref{thm:main 1} ensures that $\nabla  \bar k_n(F+\di \phi_n)=u_n$. If we denote $\bar u:= \nabla\bar k(F+\di \bar\phi)$, then, up to a subsequence, the sequence $\{u_n\}_{n\in \N}$ converges a.e to $\bar u$.
\newline
{\bf Step 2.} Let $l\in C_b(\Rd)$. The strong convergence in $L^2(\Omega)$ of $\left\lbrace\beta_n\right\rbrace_{n=1}^{\infty} $ to 1 established in Lemma \ref{lem:various compactness}  and the almost everywhere convergence of  $\{u_n\}_{n\in \N}$   to $\bar u$ obtained in Step 1 ensure that $\lim_{n\rightarrow\infty}\intom \beta_n l(u_n)dx=\intom l(\bar u (x))dx$. As $(u_n,\beta_n)\in \mathscr U_S^*$, $\intom \beta_n(x) l(u_n)dx=\intom l(y)dy$ for all $l\in C_b(\Rd)$. It follows that in the limit $\intom l(\bar u)dx=\intom l(y)dy$ for all $l\in C_b(\Rd)$  and thus $\bar u\in \mathscr U_S^1$.
\\
{\bf Step 3.} We recall that 
$$ I_{  n}(u,\beta) =V_\Set^f(u)+\intom\left(  H_{  n}(\beta)-u\cdot F\right)\;dx. $$
Since $u\mapsto V_S^f(u)$ is lower-semicontinuous as a supremun
of affine functions, by applying the Fatou's Lemma, we have
\[
\liminf_n I_{  n}(u_n,\beta_n)
\geq V_S^f(\bar u)+\intom -\bar u\cdot F\;dx =I_0(\bar u). 
\]
Let $u_0$ be the unique minimizer of $I_0$ over $\mathscr U_S^1$ as given by Theorem \ref{thm:main 2}. Then,
 \begin{equation}\label{linf In}
 \liminf_n I_{  n}(u_n,\beta_n)
 \geq  I_0(\bar u)\geq I_0(u_0).
 \end{equation}
Meanwhile, as $C_n\p C_0$ and $(k_0,l_0,\phi_0)$ is a minimizer of $J$ over  $ C_0$, we have $$J(k_0,l_0,\phi_0)\leq J(k_n,l_n,\phi_n).$$
This, along with the duality established in Theorem \ref{thm:main 1}  imply that
\begin{equation}
\begin{split}\label{eq: lim sup In}
 \limsup_n I_{ n}(u_n,\beta_n)\leq \limsup_n \left(-J(k_n,l_n,\phi_n)\right) \leq & -J(k_0,l_0,\phi_0)=I_0(u_0).
 \end{split}
\end{equation}
We combine \eqref{linf In} and \eqref{eq: lim sup In} to obtain $I_0(\bar u)=I_0(u_0)$. As $u_0$ is the unique minimizer of $I_0$  over  
 $\mathscr U_S^1$  we have $u_0=\bar u$. We note that the limit $\bar u$ does not depend on the   subsequence of $\{u_n\}_n$ chosen. Thus, the whole sequence  $\{u_n\}_n$ converges a.e. to $u_0$.
In addition, $\{I_n(u_n,\beta_n)\}_n$
converges  
to $I_0(u_0)$.
\cqfd 
\section{Appendix}
\subsection{Proof of Lemma \ref{lem:lcb}}\label{sec:app1}
We will prove Lemma \ref{lem:lcb} through 
two  lemmas. The 
results of the first lemma can be found in Lemma 4.3 of
\cite{awi-gan}. We give here a sketch of the proof for the
convenience of the reader.

\begin{lemme}\label{lem:lcb:lem1}
Assume $\mathbf{(A2)}$ holds.
Consider a lower semicontinuous  function 
$l:\mathbb{R}^d  \to \bar{\mathbb{R}}  $
such that  $\inf_{\bar\Lambda}l>-\infty$; $l$ is finite on  $\Lambda$ and $l \equiv +\infty$ on $\mathbb{R}^d  \setminus\bar\Lambda$. Set $k= l^{\#}$ and let $w\in \mathbb{R}^d  $.  Then: 
\begin{enumerate}
 \item There exist $\bar u\in\bar\Lambda$ and $\bar t>0$  such that  
\begin{equation}\label{eq:app lem 1}
 k( w ) =-\bar tl(\bar u)-H(\bar t)-\bar u\cdot w.
\end{equation}
Moreover, $ \bar u $ and $\bar t $ satisfy   $\bar u\in\pd k( w ) $ and $H'(\bar t)+l(\bar u)=0$.
\item If
$ k$ is differentiable at $w$ then $\bar u$ and $\bar t$ are uniquely determined by 
$\bar u=\nabla k(w)$ and $\bar t=(H')^{-1}(-l(\bar u))$.
\end{enumerate}
\end{lemme}
Proof. 1) We have 
\be\label{eq:app k(w)}
k(w)=\sup\{u\cdot w-l(u)t-H(t):u\in \bar\Lambda,t>0\}.
\ee
Consider a maximizing sequence
 $ \{(u_n,t_n)\}_{n=1}^{\infty} $ in \eqref{eq:app k(w)}. As $0\in \Lambda$,
we may assume without loss of generality that 
 \begin{align*}
 u_nw-l(u_n)t_n-H(t_n)&\geq 0\cdot w-l(0)-H(1)
 \end{align*}
for $n\geq 1$. It follows that
\begin{align*}
|w|r^*+l(0)+H(1)&\geq (\inf_{\bar \Lambda }l)t_n+H(t_n)
\end{align*}
for $n\geq 1$. In light of the growth condition on $ H $ in \textbf{(A2)} there exists positive real numbers $ \alpha$  such that $ \{t_n\}_{n=1}^{\infty}\subset [\alpha,\alpha^{-1}] $. As $ \Lambda $ is bounded, we may assume without loss of generality that the sequence  
 $ \{(u_n,t_n)\}_{n=1}^{\infty}$ converges to some $ (\bar u,\bar t)\in \bar\Lambda\times[\alpha,\alpha^{-1}] $.
 We next use the lower semicontinuity of $ H $ and $ l $ to deduce that
   \begin{equation}\label{eq: equilibrium}
  k(w)= \bar u\cdot w-l(\bar u)\bar t-H(\bar t).
 \end{equation}
 Note that $k(w)\geq   \bar u\cdot w-l(\bar u)t-H( t) \text{ for all }t>0$. In view of \eqref{eq: equilibrium}, it follows that  $g: (0,\infty)\to \mathbb{R}$ defined by $g(t)=\bar u\cdot w-l(\bar u)t-H( t) $ admits a maximum at $ \bar t $. As $g$ is differentiable at $\bar t$, we have $ g'(\bar t)=0 $, that is,  $ l(\bar u)+H'(\bar t)=0 $. Next, observe that $k(z)\geq   \bar u\cdot z-l(\bar u)\bar t-H( \bar t) \text{ for all }z\in \R^d$. In light of the convexity of $k$ we have that $ \bar u\in \partial k(w) $.
 
2) Assume that $k$ is differentiable at $w$. Then, $\bar u$ is uniquely determined as $ \bar u=\nabla k (w) $. As $ H'(\bar t)=-l_0(\bar u) $ and $ H' $ is a bijection, we obtain that $ \bar t $ is also uniquely determined as $\bar t=(H')^{-1}(-l(\bar u))$.
\qed

The second lemma which is inspired by Lemma 4.4 in \cite{awi-gan} is the following:

\begin{lemme}\label{lem:lcb lem2}
Assume assumption $\mathbf{(A2)}$ holds.
Consider a lower semicontinuous  function 
$l_0:\mathbb{R}^d  \to \bar{\mathbb{R}}  $
   such that  $\inf_{\bar\Lambda}l_0>-\infty$; $l_0$ is finite on  $\Lambda$ and $l_0\equiv +\infty$ on $\mathbb{R}^d  \setminus\bar\Lambda$. Set $k_0= ({l_0})^{\#}$.
   Let $\hat l\in C_b(\mathbb{R}^d  )$ and let $1\geq \epsilon>0$. Define $l_\epsilon=l_0+\epsilon \hat l$ and $k_\epsilon={\left(l_\epsilon  \right)}^{\#}$. Let $v\in \mathbb{R}^d  $  be such that    $ k_0$ is differentiable at $v$. 
 \begin{enumerate}
  \item \label{lcbiia:app}There exists a constant $M$ independent of $v$ and $\epsilon$ 
  such that
\begin{equation}\label{eq:ratio bound 0}
  \left|\frac{k_\epsilon(v)-k_0( v ) }{\epsilon} \right|\leq M.
 \end{equation}

  \item \label{lcbiib:app}We have 
  \begin{equation}\label{eq:limit ratio bound}
  \lim _{\epsilon\rightarrow 0 }   \frac{k_\epsilon(v)-k_0( v ) }{\epsilon}=- t_0\hat l(u_0).
  \end{equation}
 \end{enumerate} 

\end{lemme} 
Proof. Note that the map $ l_\epsilon=l_0+\epsilon \hat l $  is bounded below by $ m-|\hat l|_\infty $. As $k_\epsilon={\left(l_\epsilon  \right)}^{\#}$ and  $k_0={\left(l_0  \right)}^{\#}$, lemma \ref{lem:lcb:lem1} ensures that there exist $t_0, t_\epsilon>0 $ and $ u_0, u_\epsilon\in \bar \Lambda $ such that 
$$k_\epsilon(v)=u_\epsilon v-l(u_\epsilon)t_\epsilon-H(t_\epsilon)$$
and
$$k_0(v)=u_0 v-l(u_0)t_0-H(t_0).$$
We then have
\begin{equation}\label{eq: equili1}
 k_\epsilon(v) =-\epsilon \hat l(u_\epsilon)t_\epsilon
+u_\epsilon v-l_0(u_\epsilon)t_\epsilon-H(t_\epsilon)\leq -\epsilon\hat  l(u_\epsilon)t_\epsilon +k_0(v)
\end{equation}
and
\begin{equation}\label{eq: equili2}
 k_0(v)=\epsilon \hat l(u_0)t_0
+u_0 v-l_\epsilon(u_0)t_0-H(t_0)\leq \epsilon \hat l(u_0)t_0 +k_\epsilon(v).
\end{equation}
 We combine \eqref{eq: equili1} and \eqref{eq: equili2} to get 
\begin{equation}\label{eq: ratio bound}
-\hat l(u_0)t_0\leq   \frac{(k_\epsilon (v)- k_0(v))}{\epsilon}\leq 
-\hat l(u_\epsilon)t_\epsilon.
\end{equation}
Using again lemma \ref{lem:lcb:lem1} we have 
$$t_{\delta}= (H')^{-1}( l_{\delta}(u_{\delta})),\quad u_{\delta}\in \partial k_{\delta}(v)\qquad \delta\in \{0,\epsilon\}.$$
As $l_{\delta}$ is bounded below by $m-|l|_{\infty}$,  we use the fact that $ H' $ is a continuous and strictly increasing bijection from $(0,\infty)$ to $\mathbb R$ to deduce that $ t_\delta $ is bounded above by $ M_1>0 $ given by $M_1:=(H')^{-1}(-m+|\hat l|_\infty)$. This bound on $t_{\delta}$ combined with \eqref{eq: ratio bound} yields a constant  $M:=|\hat l|_\infty(H')^{-1}(-m+|\hat l|_\infty)$ such that \eqref{eq:ratio bound 0} holds. As a result
$\lim_{\epsilon\to 0^+}k_\epsilon(v)=k_0(v).$
Next, let $ \{e_n\}_{n=1}^{\infty}\subset (0,1] $ converging to $ 0 $ such that $ \limsup_{\epsilon\to 0} \hat l(u_\epsilon)t_\epsilon =\lim_{n\to \infty}\hat l(u _{e_n}
)t_{e_n}$. Without loss of generality, we may assume that $ \{u _{e_n}\}_{n=1}^{\infty}
 $ converges to some $ \bar u\in \bar \Lambda $ and $\{ t _{e_n}\}_{n=1}^{\infty} 
 $ converges to $ \bar t \in [0,M_1]$.
Exploiting the lower semicontinuity of $ l_0 $, $ \hat l $ and $ H $, we get:
\begin{align*}
k_0(v)&=\lim_{n\to \infty} k_{e_n}(v)
\\
&=\lim_{n\to \infty} u_{e_n}v-l_{e_n}(u_{e_n})t_{e_n}-H(t_{e_n})
\\
&\leq \bar u v-l_0(\bar u)\bar t-H(\bar t)
\\
&\leq k_0(v).
\end{align*}
It follows that $ k_0(v)=\bar u v-l_0(\bar u)\bar t-H(\bar t) $. As $k_0$ is differentiable at $v$, we have $ t_0=\bar t $ and $ u_0=\bar u $.
We use \eqref{eq: ratio bound}, the definition of $ \{e_n\}_{n=1}^{\infty}$, the convergence of $ \{u _{e_n}\}_{n=1}^{\infty}
 $ and $\{ t _{e_n}\}_{n=1}^{\infty} $ to obtain 
\begin{equation}\label{eq: limit tl(u)}
-t_0\hat l(u_0)\leq \liminf_{\epsilon \to 0}-t_\epsilon \hat l(u_\epsilon)\leq \limsup_{\epsilon \to 0}-t_\epsilon \hat l(u_\epsilon)=\lim_{n\to \infty}-t_{e_n}\hat l(u_{e_n})=-t_0\hat l(u_0).
\end{equation}
As a result, $\lim_{\epsilon \to 0}-t_\epsilon \hat l(u_\epsilon)=-t_0\hat l(u_0)  $. We invoke one more time equation \eqref{eq: ratio bound} to obtain \eqref{eq:limit ratio bound}.

\subsection{Some properties of the Legendre transform.}
We have the following Lemma which is similar to
Lemma  \ref{lem:lcb} but uses the Legendre transform instead of the $(\cdot)^\#$ operator.
\begin{lemme}\label{lem:lcb:H0}
Consider a lower semicontinuous  function 
$l_0:\mathbb{R}^d  \to \bar{\mathbb{R}}  $
   such that  $\inf_{\bar\Lambda}l_0>-\infty$; $l_0$ is finite on  $\Lambda$ and $l_0\equiv +\infty$ on $\mathbb{R}^d  \setminus\bar\Lambda$. Set $k_0= ({l_0})^*$. 
\begin{enumerate}
\item There exists a measurable map $ T_0:\mathbb R^d \to \mathbb R^d$ such that 
$ k_0(v)=v\cdot T_0(v)-l_0(T_0(v)) $
for all $ v\in \mathbb R^d $
and $ T_0(v)=\nabla k_0(v) $ whenever $ k_0 $ is differentiable at $ v\in \mathbb R^d $.
 \item \label{lcbii:H0} Let $\hat l\in C_b(\mathbb{R}^d  )$ and let $1\geq \epsilon>0$. Define $l_\epsilon=l_0+\epsilon \hat l$ and $k_\epsilon={\left(l_\epsilon  \right)^*}$.
 \begin{enumerate}
 
  \item \label{lcbiia:H0}
  For all $ v\in \mathbb R^d $ we have :
 \[
  \left|\frac{k_\epsilon(v)-k_0( v ) }{\epsilon}\right|\leq |\hat l|_\infty.
 \]
  \item \label{lcbiib_H0}
  For $\epsilon\in(0,1)$, there exists a  
  map $ T_\epsilon:\mathbb R^d \to \mathbb R^d $ satisfying 
for all $ v\in \mathbb R^d $:
$
k_\epsilon(v)=vT_\epsilon(v)-l_\epsilon(T_\epsilon(v))$.
When $ k_0 $ is differentiable at $ v\in \mathbb R^d $, we have
  $ \lim_{\epsilon \to 0}T_\epsilon(v)= \nabla k_0(v)$ and
  \[ \lim _{\epsilon\rightarrow 0 }   \frac{k_\epsilon(v)-k_0( v ) }{\epsilon}=- t_0\hat l(\nabla k_0(v)).\]
 \end{enumerate}

\end{enumerate}

\end{lemme} 
Proof. 1.) Let $ v\in \mathbb{R}^d $.
We have
\[k_0(v)=\sup\{uv-l_0(u): u\in \mathbb R^d\}=\sup\{uv-l_0(u): u\in \bar\Lambda\}.\]
We use the lower semicontinuity of $ l_0 $ and the compactness of $ \bar \Lambda $ to deduce that there exists $ \bar u\in \Lambda $ such that $ k_0(v)=\bar u  v-l_0(\bar u). $
We have $ k_0(w)-(\bar uw-l_0(\bar u))\geq 0 $ for all $ w\in \mathbb R^d $ while $ k_0(v)-(\bar uv-l_0(\bar u))= 0 $.
Since 
$ k_0 $ is convex, we deduce that $ \bar u \in \partial k_0(v) $. 
\\
Next, for
$ v\in \mathbb{R}^d $, define
\[\Gamma(v)=\{u\in \bar\Lambda :k_0(v)=uv-l_0(u) \}.\]
Assume $ \{u_n\}_{n\in \mathbb N}\subset \mathbb{R}^d $ converges to $ u $; $ \{v_n\}_{n\in \mathbb N}\subset \mathbb{R}^d $ converges to $ v $ and for all $ n\in \mathbb N $, one has $ u_n\in \Gamma(v_n) $. Then $ u\in \Gamma(v) $. Indeed, one has
\[k_0(v)\leq\liminf_{n\to \infty}k_0(v_n)=\liminf_{n\to \infty}\left(u_nv_n-l_0(u_n)\right)\leq uv-l_0(u)\leq k_0(v). \]
Therefore, $ uv-l_0(u)=k_0(v) $
and $ u\in \Gamma(v) $. As a result, the multifunction $ \Gamma:\mathbb R^d \rightrightarrows \mathbb R^d $ is closed and nonempty valued.
By the Measurable Selection Theorem [Corollary 14.6, \cite{rockafellar2009variational}], there exists a measurable map $ T_0:\mathbb R^d\to \mathbb R^d $ such that for all $ v\in \mathbb R^d $, one has $ T_0(v)\in \Gamma(v) $. That is $ k_0(v)=vT_0(v)-l_0(T_0(v)) $. As $ T(v)\in \Gamma(v)\subset \partial k_0(v) $, we also have $ T_0=\nabla k_0 $ almost everywhere. 
\\
2.) For $ \epsilon>0 $, $ l_\epsilon $ is bounded below
and satisfies the  
hypothesis on $ l_0 $. Let 
$ k_\epsilon=l_\epsilon ^* $ and consider a map $ T_\epsilon $ satisfying 
for all $ v\in \mathbb R^d $:
$
k_\epsilon(v)=vT_\epsilon(v)-l_\epsilon(T_\epsilon(v))$
as given by part 1.).
We have for $ v\in \mathbb R^d $:
\begin{equation}\label{eq:kepsilon ineq0}
k_\epsilon(v)=vT_\epsilon(v)-l_\epsilon(T_\epsilon(v))=-\epsilon \hat l(T_\epsilon(v))+vT_\epsilon(v)-l_0(T_\epsilon(v))
\leq -\epsilon \hat l(T_\epsilon(v))+k_0(v).
\end{equation}
 Similarly, for $ v\in \mathbb R^d $ we have
\begin{equation} \label{eq:kepsilon ineq1}
k_0(v)=vT_0(v)-l_0(T_0(v))=\epsilon \hat l(T_0(v))+vT_0(v)-l_\epsilon(T_0(v))\leq \epsilon \hat l(T_0(v))+k_\epsilon(v).
\end{equation}
We combine \eqref{eq:kepsilon ineq0} and \eqref{eq:kepsilon ineq1} to get
\begin{equation}\label{eq:kepsilon:H0}
-\hat l(T_0(v))\leq
\frac{k_\epsilon(v)-k_0(v)}{\epsilon}\leq -\hat l(T_\epsilon(v)),
\end{equation}
which leads to
\begin{equation}\label{eq:kepsilonBound:H0}
\left|\frac{k_\epsilon(v)-k_0(v)}{\epsilon}\right|\leq |\hat l|_\infty.
\end{equation}
Consider a sequence
$ \{\epsilon_n\}_n $ converging to $ 0 $. The sequence 
$ \{T_{\epsilon_n}(v)\}_n $ is bounded so we may find a subsequence $ \{\epsilon'_n\}_n $ of $ \{\epsilon_n\}_n $ such that the sequence $ \{T_{\epsilon'_n}(v)\}_n $ converges to $ u\in \bar \Lambda $. We then have:
\begin{equation}\label{eq:maximizer  sup k_0}
k_0(v)=\lim_{n\to \infty}k_{\epsilon_n'}(v)=
\lim_{n\to \infty}\left(vT_{\epsilon_n'}(v)-l_{\epsilon_n}(T_{\epsilon_n'}(v))\right)
\leq vu-l_0(u)\leq k_0(v).
\end{equation}
We use  \eqref{eq:maximizer  sup k_0} to obtain $ k_0(v)=vu-l_0(u)$ and thus $ u=\nabla k_0(v) $ as $ k_0 $ is differentiable at $ v $.
It follows that $ \lim_{\epsilon \to 0}T_\epsilon(v)= \nabla k_0(v)$.
We use Equation \eqref{eq:kepsilon:H0} and the continuity of $ \hat l $ to obtain
\[ \lim _{\epsilon\rightarrow 0 }   \frac{k_\epsilon(v)-k_0( v ) }{\epsilon}=- t_0 \hat l(\nabla k_0(v)).\]
\qed


\bibliographystyle{spmpsci}

\begin{thebibliography}{10}

\bibitem{ambrosio2005gradient}
L.~Ambrosio, N.~Gigli, and G.~Savar{\'e}.
\newblock {\em Gradient Flows: In Metric Spaces And In The Space Of Probability
  Measures}.
\newblock Lectures in mathematics. Birkh{\"a}user, 2005.
\bibitem{awi-gan}
R.~Awi and W.~Gangbo.
\newblock A polyconvex integrand; euler–lagrange equations and uniqueness of
  equilibrium.
\newblock {\em Archive for Rational Mechanics and Analysis}, 214(1):143--182,
  2014.

\bibitem{Brenier1991}
Y.~{Brenier}.
\newblock Polar factorization and monotone rearrangement of vector-valued
  functions.
\newblock {\em Comm. Pure Appl. Math.}, 44:375--417, 1991.

\bibitem{Dac_dmcv}
B.~Dacorogna.
\newblock {\em Direct methods in the Calculus of Variations}.
\newblock Springer, 2008.

\bibitem{egs}
L.~C. Evans, W.~Gangbo, and O.~Savin.
\newblock Diffeomorphisms and nonlinear heat flows.
\newblock {\em SIAM J. Math. Anal.}, 2005.

\bibitem{Evans20120339}
L.~C.~Evans.
\newblock Monotonicity formulae for variational problems.
\newblock {\em Philosophical Transactions of the Royal Society of London A:
  Mathematical, Physical and Engineering Sciences}, 371(2005), 2013.

\bibitem{mtfpf}
L.~C.~Evans and R.~Gariepy.
\newblock {\em Measure Theory and Fine Properties of Functions}.
\newblock CRC Press, 1992.

\bibitem{Gangbo1994}
W.~{Gangbo}.
\newblock An elementary proof of the polar factorization of vector-valued
  functions.
\newblock {\em Archive for Rational Mechanics and Analysis}, 128:381--399,
  1994.

\bibitem{gangbomccann}
W.~Gangbo and R.J.~McCann.
\newblock The geometry of optimal transportation.
\newblock {\em Acta Math.}, 1996.

\bibitem{Gan_uecsc}
W.~Gangbo and R.~Van~Der Putten.
\newblock Uniqueness of equilibrium configurations in solid crystals.
\newblock {\em SIAM J. MATH. ANAL.}, 2000.

\bibitem{gan-ma}
W.~Gangbo.
\newblock The monge mass transfer problems and its applications.
\newblock {\em NSF-CBMS Conference. Contemporary Mathematics}, 1999.

\bibitem{gan-wes}
W.~Gangbo and M.~Westdickenberg.
\newblock {Optimal Transport for the System of Isentropic Euler Equations}.
\newblock {\em Communications in Partial Differential Equations},
  34:1041--1073, 2009.
\bibitem{Roc_ca}
R.~T. Rockafellar.
\newblock {\em Convex Analysis}.
\newblock Princeton Univ. Press, 1970.

\bibitem{villani2003topics}
C.~Villani.
\newblock {\em Topics in Optimal Transportation}.
\newblock Graduate Studies in Mathematics Series. Amer Mathematical Society,
  2003.
  
 \bibitem{rockafellar2009variational}
~Rockafellar, R.T. and ~Wets, M. and ~Wets, R.J.B.
\newblock {\em Variational Analysis}.
\newblock Grundlehren der mathematischen Wissenschaften,
  2009.

\end{thebibliography}

\end{document}